\documentclass[a4paper,11pt]{article}
\usepackage{}
\usepackage{mathrsfs}
\usepackage{amsfonts}
\usepackage{amsfonts}
\usepackage{amssymb}
\usepackage{mathrsfs}
\usepackage{indentfirst,latexsym,bm,amsmath,amssymb,amsthm}
\usepackage{xcolor}

\usepackage{indentfirst,latexsym,bm,amsmath,amssymb,amsthm}
\usepackage[dvips]{graphicx}

\makeatletter\@addtoreset{equation}{section} \makeatother
\newtheorem{thm}{Theorem}[section]
\newtheorem{Lemma}{Lemma}[section]

\makeatletter 

\allowdisplaybreaks
\makeatother
\title{Global existence for initial-boundary value problems of\\ one-dimension quasilinear wave equations with null conditions}

\author{Dongbing Zha\thanks{  Department of Mathematics, Donghua University, Shanghai 201620, PR China.{ E-mail address: ZhaDongbing@163.com}}}

\begin{document}

\maketitle
\begin{abstract}
We consider the initial-boundary value problems on $\mathbb{R}^{+}\times \mathbb{R}^{+}$ for one-dimension systems of quasilinear wave equations with null conditions. We show that for homogeneous Dirichlet boundary values and sufficiently small initial data, classical solutions always globally exist. The key innovation in the proof is a new framework of bootstrap argument via coupled high-low order energy estimates.
\\
\emph{keywords}: One-dimension quasilinear wave equations; initial-boundary value problem; null condition; global existence.\\
\emph{2020 MSC}: 35L05; 35L53; 35L72.
\end{abstract}

\pagestyle{plain} \pagenumbering{arabic}

\section{Introduction }

For the Cauchy problem of one-dimension semilinear wave equations with null conditions, \cite{MR3783412} shows the global existence of classical solutions with small initial data.  This result strengthens a previous one in \cite{MR3121700}, which shows the global existence under
some stronger null conditions. One-dimension nonlinear wave equation is essentially a kind of transport system, which does not admit decay in time. The mechanism for the global existence in such case is the interaction of waves with different speeds, which will lead to the decay in time of nonlinear terms. In order to display this mechanism,  a kind of weighted energy estimates with positive weights is developed in \cite{MR3783412}.
This mechanism is completely different from the corresponding ones in high dimensional cases. We refer the reader to
Christodoulou and Klainerman's pioneering works  for the global existence of classical solutions for nonlinear wave equations with null conditions in three space dimensions \cite{MR820070, Klainerman86}, and of Alinhac's global existence result for the case of two space dimensions \cite{Alinhac01} (see also \cite{MR3729247,MR3912654} for some thorough studies), which are all based on the decay in time of liner waves.

Then for the Cauchy problem of one-dimension quasilinear wave equations with null conditions, the global existence of classical solutions in the small data setting is shown in \cite{MR4098041}. The proof in \cite{MR4098041} is based on weighted energy estimates with positive weights in \cite{MR3783412},  some space-time weighted energies and new observations concerning the null structure in the quasilinear part.

%
Based on the above global existence results in the Cauchy problem case for one-dimension nonlinear wave equations, in this paper, we intend to consider the corresponding topic in the initial-boundary value problem case.

 The outline of this paper is as follows. The remainder of this introduction will be devoted to the description of a statement
of main result. In Section \ref{hjsedddd}, some necessary tools used to prove Theorem \ref{mainthm} are introduced.
Section \ref{xhzuyo898} is devoted to the proof of Theorem \ref{mainthm}.

\subsection{Main result}
Let $(t,x)$ denote the usual Cartesian coordinates in $\mathbb{R}^{+}\times\mathbb{R}^{+}$, and define also the null coordinates
\begin{align}
\xi=\frac{t+x}{2},~~\eta=\frac{t-x}{2},
\end{align}
as well as the corresponding null vector fields
\begin{align}\label{xuyao890hj}
\partial_{\xi}=\partial_t+\partial_x,~~\partial_{\eta}=\partial_t-\partial_x.
\end{align}
We also denote briefly $u_{\xi}=\partial_{\xi}u$ and $u_{\eta}=\partial_{\eta}u$.\par
Consider the following one-dimension system
\begin{align}\label{quasiwave}
u_{\xi\eta}&=A_1(u,u_{\xi},u_{\eta})u_{\xi\eta}
+A_2(u,u_{\xi},u_{\eta})u_{\xi\xi}+A_3(u,u_{\xi},u_{\eta})u_{\eta\eta}+F(u,u_{\xi},u_{\eta}),
\end{align}
where the unknown function $u=u(t,x): \mathbb{R}^{+}\times\mathbb{R}^{+}\longrightarrow \mathbb{R}^{n}$, for $i=1,2,3$, $A_i: \mathbb{R}^{n}\times \mathbb{R}^{n}\times \mathbb{R}^{n}\longrightarrow  \mathbb{R}^{n\times n}$ are given smooth and matrix valued functions, and $F: \mathbb{R}^{n}\times \mathbb{R}^{n}\times \mathbb{R}^{n}\longrightarrow  \mathbb{R}^{n}$ is a given smooth and vector valued function. Moreover, we will always assume that $A_i ~(i=1,2,3)$ are symmetric, which means that they take values of symmetric matrixes.

Note that the left-hand side $u_{ \xi \eta } = u_{ t t } - u_{ x x }$ of \eqref{quasiwave} represents the standard one-dimension wave operator, while the right-hand side of \eqref{quasiwave} represents nonlinear perturbations.
Any one-dimension systems of quasilinear wave equations admit the form of \eqref{quasiwave}. We note that in one space dimension case waves
do not decay, and any nonlinear resonance (even arbitrarily high order) can lead to finite time
blowup. Thus in order to avoid the formation of singularities (in the small data setting), the null conditions
should be imposed on all order nonlinearities. Klainerman first suggested such structural condition
for three dimensional nonlinear wave equations \cite{Klainerman82}. It means that any given plane wave solution of
the homogeneous linear wave equation (i.e., left traveling wave $f(\xi)$ and right traveling wave
$g(\eta)$ in one-dimension case) also satisfies the nonlinear wave equation.

We call that the system \eqref{quasiwave} satisfies the null conditions, according to the above definition, if near the origin in $\mathbb{R}^{n}\times \mathbb{R}^{n}\times \mathbb{R}^{n}$, it holds that
\begin{align}\label{order1}
A_1(u,u_{\xi},u_{\eta})&=\mathscr{O}(|u|+|u_{\xi}|+|u_{\eta}|),\\\label{order22}
A_2(u,u_{\xi},u_{\eta})&=\mathscr{O}(|u_{\eta}|),\\\label{order33}
A_3(u,u_{\xi},u_{\eta})&=\mathscr{O}(|u_{\xi}|),\\\label{order4}
F(u,u_{\xi},u_{\eta})&=\mathscr{O}(|u_{\xi}||u_{\eta}|).
\end{align}

We point out that many important equations arising from geometry
and physics in one space dimension, such as timelike minimal surface \cite{MR2045426}, Faddeev model \cite{MR3585834},  Chaplygin gas \cite{MR694243}, Born-Infeld model \cite{10.2307/2935568}, wave maps \cite{MR596432}, fall into the form of \eqref{quasiwave} with structural conditions \eqref{order1}--\eqref{order4}.

 The purpose of this paper is to treat the initial-boundary value problem on $\mathbb{R}^{+}\times \mathbb{R}^{+}$ for one-dimension system of quasilinear wave equations
\eqref{quasiwave} under the null conditions.

Consider the homogeneous Dirichlet boundary condition
\begin{align}\label{boundary}
u(t,0)=0,~t\geq 0,
\end{align}
and the initial condition
\begin{align}\label{initial}
t=0: u=u_0(x),~ u_{t}=u_{1}(x), x\geq 0.
\end{align}
As usual, we also always assume that the initial data satisfy the compatibility conditions of order three. That means, $u_0(0)=0, u_1(0)=0$, $u'_0(0)=0$, $u'_1(0)=0$, $u^{''}_0(0)=0$ and $u^{''}_1(0)=0$.

The main result of this paper is the following
\begin{thm}\label{mainthm}
For the system \eqref{quasiwave}, assume that $A_1, A_2, A_3$ are symmetric, \eqref{order1}, \eqref{order22}, \eqref{order33}, \eqref{order4} hold.
Then for all $0<\delta<1$, there exist a positive constant $\varepsilon_0$ such that for any $0<\varepsilon\leq \varepsilon_0$, if
\begin{align}\label{xxxjkdddd900}
\sum_{l=0}^{4}\|\langle x\rangle^{3+3\delta}\partial_x^{l}u_0\|_{L_{x}^2(\mathbb{R}^{+})}+\sum_{l=0}^{3}\|\langle x\rangle^{3+3\delta}\partial_x^{l}u_1\|_{L_{x}^2(\mathbb{R}^{+})}\leq \varepsilon,
\end{align}
then the initial-boundary value problem \eqref{quasiwave}, \eqref{boundary}, \eqref{initial} admits a unique global classical solution.
\end{thm}


The semilinear version of Theorem \ref{mainthm} is verified in \cite{MR4191336}. The purpose of this paper is to investigate the quasilinear case, which is much more complicated than the semilinear one.
 Compared with the semilinear case, the main difficulty in the quasilinear case lies in that there will be some uncontrollable terms in the quasilinear part after integrating by parts argument, under the null conditions \eqref{order22} and \eqref{order33}.\footnote{Under some stronger null conditions, this difficulty can be avoided. See \cite{MR4612357}.} In the Cauchy problem case \cite{MR4098041}, we can solve this problem by using the null coordinates $Z=(\partial_{\xi},\partial_{\eta})$ as the commuting vector field, and treating such terms by the equation itself (integrating by parts is not used). But in the initial-boundary value problem case, this approach obviously does not work. Because we can only use the time derivative, which preserves the homogeneous boundary conditions, as the commuting
vector field. Thus in order to prove Theorem \ref{mainthm}, we will develop some new framework, the key innovation in which is a coupled high-low order energy estimate argument.

We also point out that the method of proof for Theorem {\rm{\ref{mainthm}}} can be directly used to treat the initial-boundary value problem outside of a bounded internal, which is also called the exterior domain problem. In the 3-D case, for the exterior domain problem, i.e., initial-boundary value problem outside of some compact obstacle, the corresponding analogue of Christodoulou and Klainerman's global existence results for nonlinear wave equations with null conditions {\rm{\cite{MR820070, Klainerman86}}} are obtained in {\rm{\cite{MR1887632,MR3765760,MR2110542,MR2142333}}}, etc. While, in the 2-D case, how to get the exterior domain problem analogue of Alinhac's global existence result {\rm{\cite{Alinhac01}}} is still open until now.

\section{Preliminaries}\label{hjsedddd}
\subsection{Notations}\label{hjsedddffffffd}
For the convenience, we first introduce some notations.
Fix $0<\delta<1$. Denote
\begin{align}\label{xuyao00099}
\phi(x)=\langle x\rangle^{2+2\delta},~~\theta(x)=\langle x\rangle^{6+6\delta}.
\end{align}
It is easy to verify that
\begin{align}\label{20p06}
|\phi'(x)|\leq 4\langle x\rangle^{1+2\delta},~~|\theta'(x)|\leq 12\langle x\rangle^{5+6\delta}.
\end{align}
Denote
\begin{align}
\psi(x)=\exp({-\int_{-\infty}^{x}{\langle \rho\rangle^{-(1+\delta)}}d\rho}).
\end{align}
We can verify that
\begin{align}
\psi'(x)=-\psi(x)\langle x\rangle^{-(1+\delta)}.
\end{align}
We note that there exists a positive constant $c$ such that
\begin{align}\label{rg56}
c^{-1}\leq \psi(x)\leq c,
\end{align}
thus
\begin{align}\label{here456}
c^{-1}\langle x\rangle^{-(1+\delta)}\leq -\psi'(x)\leq c\langle x\rangle^{-(1+\delta)}.
\end{align}

Now we will introduce some notations for energies.
Denote the vector field $Z=(\partial_{\xi},\partial_{\eta})$. For multi-index $a=(a_1,a_2)$, denote $Z^{a}=\partial_{\xi}^{a_1}\partial_{\eta}^{a_2}$ and $|a|=a_1+a_2$.
We will use the following weighted energies with positive weights
\begin{align}\label{Lu}
&E(u(t))=\|\langle \xi\rangle^{1+\delta}u_{\xi}\|_{L_{x}^2(\mathbb{R}^{+})}^2+\|\langle \eta\rangle^{1+\delta}u_{\eta}\|_{L_{x}^2(\mathbb{R}^{+})}^2,\\\label{Lu2}
&\mathbb{E}(u(t))=\|\langle \xi\rangle^{3+3\delta}u_{\xi}\|_{L_{x}^2(\mathbb{R}^{+})}^2+\|\langle \eta\rangle^{3+3\delta}u_{\eta}\|_{L_{x}^2(\mathbb{R}^{+})}^2,
\end{align}
and corresponding higher order energies
\begin{align}
E_k(u(t))=\sum_{|a|\leq k-1}E(Z^{a}u(t)),~~\mathbb{E}_k(u(t))=\sum_{|a|\leq k-1}\mathbb{E}(Z^{a}u(t)).
\end{align}
Denoting $S_t=[0,t]\times \mathbb{R}^+$, we further introduce the following space-time weighted energies
\begin{align}\label{Lu3}
&\mathcal{E}(u(t))=\|\langle \eta\rangle^{-\frac{1+\delta}{2}}\langle \xi\rangle^{1+\delta}u_{\xi}\|^2_{L^2_{s,x}(S_t)}+\|\langle \xi\rangle^{-\frac{1+\delta}{2}}\langle \eta\rangle^{1+\delta}u_{\eta}\|^2_{L^2_{s,x}(S_t)},\\\label{Lu4}
&\mathscr{E}(u(t))=\|\langle \eta\rangle^{-\frac{1+\delta}{2}}\langle \xi\rangle^{3+3\delta}u_{\xi}\|^2_{L^2_{s,x}(S_t)}+\|\langle \xi\rangle^{-\frac{1+\delta}{2}}\langle \eta\rangle^{3+3\delta}u_{\eta}\|^2_{L^2_{s,x}(S_t)},
\end{align}
and higher order energies
\begin{align}
\mathcal{E}_k(u(t))=\sum_{|a|\leq k-1}\mathcal{E}(Z^{a}u(t)),~~\mathscr{E}_k(u(t))=\sum_{|a|\leq k-1}\mathscr{E}(Z^{a}u(t)).
\end{align}

We point out that the weighed energy with positive weights \eqref{Lu} is first used in \cite{MR3783412}, while the space-time weighted energy \eqref{Lu3} is introduced in \cite{MR4098041}, which is inspired by \cite{Alinhac01} and \cite{MR2680391}. \eqref{Lu} and \eqref{Lu3} will be employed in our high order energy estimate. In order to display the lower order energy estimate, we will use \eqref{Lu2} and \eqref{Lu4}.

For $k\geq 2$, we can split
\begin{align}\label{GHUER444}
{E_k}(u(t))=\overline{E}_k(u(t))+\overline{\overline{{E}}}_k(u(t)),~~{\mathbb{E}_k}(u(t))=\overline{\mathbb{E}}_k(u(t))+\overline{\overline{{\mathbb{E}}}}_k(u(t)),
\end{align}
with
\begin{align}
&\overline{E}_k(u(t))=\sum_{l=0}^{k-1}\|\langle \xi\rangle^{1+\delta}\partial^l_{\xi}u_{\xi}\|_{L_{x}^2(\mathbb{R}^{+})}^2+\sum_{l=0}^{k-1}\|\langle \eta\rangle^{1+\delta}\partial^l_{\eta}u_{\eta}\|_{L_{x}^2(\mathbb{R}^{+})}^2,\\
&\overline{\mathbb{E}}_k(u(t))=\sum_{l=0}^{k-1}\|\langle \xi\rangle^{3+3\delta}\partial^l_{\xi}u_{\xi}\|_{L_{x}^2(\mathbb{R}^{+})}^2+\sum_{l=0}^{k-1}\|\langle \eta\rangle^{3+3\delta}\partial^l_{\eta}u_{\eta}\|_{L_{x}^2(\mathbb{R}^{+})}^2,
\end{align}
and
\begin{align}
\overline{\overline{{E}}}_k(u(t))
&=\sum_{|a|\leq k-2}\|\langle \xi\rangle^{1+\delta}Z^{a}u_{\xi\eta}\|_{L_{x}^2(\mathbb{R}^{+})}^2
+\sum_{|a|\leq k-2}\|\langle \eta\rangle^{1+\delta}Z^{a}u_{\xi\eta}\|_{L_{x}^2(\mathbb{R}^{+})}^2,\nonumber\\
\overline{\overline{{\mathbb{E}}}}_k(u(t))
&=\sum_{|a|\leq k-2}\|\langle \xi\rangle^{3+3\delta}Z^{a}u_{\xi\eta}\|_{L_{x}^2(\mathbb{R}^{+})}^2
+\sum_{|a|\leq k-2}\|\langle \eta\rangle^{3+3\delta}Z^{a}u_{\xi\eta}\|_{L_{x}^2(\mathbb{R}^{+})}^2.
\end{align}
Correspondingly, for the space-time weighted energy, we can also write
\begin{align}\label{GHUEFFFFR444}
{\mathcal {E}_k}(u(t))=\overline{\mathcal {E}}_k(u(t))+\overline{\overline{{\mathcal {E}}}}_k(u(t)),~~
{\mathscr{E}_k}(u(t))=\overline{\mathscr{E}}_k(u(t))+\overline{\overline{{\mathscr{E}}}}_k(u(t)),
\end{align}
with
\begin{align}
\overline{\mathcal {E}}_k(u(t))=&\sum_{l=0}^{k-1}\|\langle \eta\rangle^{-\frac{1+\delta}{2}}\langle \xi\rangle^{1+\delta}\partial^{l}_{\xi}u_{\xi}\|^2_{L^2_{s,x}(S_t)}+\sum_{l=0}^{k-1}\|\langle \xi\rangle^{-\frac{1+\delta}{2}}\langle \eta\rangle^{1+\delta}\partial_{\eta}^{l}u_{\eta}\|^2_{L^2_{s,x}(S_t)},\\
\overline{\mathscr{E}}_k(u(t))=&\sum_{l=0}^{k-1}\|\langle \eta\rangle^{-\frac{1+\delta}{2}}\langle \xi\rangle^{3+3\delta}\partial^{l}_{\xi}u_{\xi}\|^2_{L^2_{s,x}(S_t)}+\sum_{l=0}^{k-1}\|\langle \xi\rangle^{-\frac{1+\delta}{2}}\langle \eta\rangle^{3+3\delta}\partial_{\eta}^{l}u_{\eta}\|^2_{L^2_{s,x}(S_t)},
\end{align}
and
\begin{align}
&\overline{\overline{{\mathcal {E}}}}_k(u(t))\nonumber\\
&=\sum_{|a|\leq k-2}\|\langle \eta\rangle^{-\frac{1+\delta}{2}}\langle \xi\rangle^{1+\delta}Z^{a}u_{\xi\eta}\|^2_{L^2_{s,x}(S_t)}+\sum_{|a|\leq k-2}\|\langle \xi\rangle^{-\frac{1+\delta}{2}}\langle \eta\rangle^{1+\delta}Z^{a}u_{\xi\eta}\|^2_{L^2_{s,x}(S_t)},\\
&\overline{\overline{{\mathscr{E}}}}_k(u(t))\nonumber\\
&=\sum_{|a|\leq k-2}\|\langle \eta\rangle^{-\frac{1+\delta}{2}}\langle \xi\rangle^{3+3\delta}Z^{a}u_{\xi\eta}\|^2_{L^2_{s,x}(S_t)}+\sum_{|a|\leq k-2}\|\langle \xi\rangle^{-\frac{1+\delta}{2}}\langle \eta\rangle^{3+3\delta}Z^{a}u_{\xi\eta}\|^2_{L^2_{s,x}(S_t)}.
\end{align}

In the energy estimates, we can only use the time derivative as the commuting vector field. Thus we further introduce
\begin{align}
\widetilde{{E}}_k(u(t))=\sum_{l=0}^{k-1}{E}(\partial_t^{l}u(t)),~~\widetilde{{\mathbb{E}}}_k(u(t))=\sum_{l=0}^{k-1}{\mathbb{E}}(\partial_t^{l}u(t)),
\end{align}
and
\begin{align}
\widetilde{\mathcal {E}}_k(u(t))=\sum_{l=0}^{k-1}{\mathcal {E}}(\partial_t^{l}u(t)),~~\widetilde{\mathscr{E}}_k(u(t))=\sum_{l=0}^{k-1}{\mathscr{E}}(\partial_t^{l}u(t)).
\end{align}

\subsection{Pointwise estimates}
The following two lemmas follow from the fundamental theorem of calculus, chain rule and Leibniz's rule.
\begin{Lemma}\label{xuyaouioo999}
 Assume that $u$ is smooth, $A_1=A_1(u,u_{\xi},u_{\eta}), A_2=A_2(u,u_{\xi},u_{\eta})$, $A_3=A_3(u,u_{\xi},u_{\eta})$ satisfies \eqref{order1}, \eqref{order22},
\eqref{order33}, respectively, and
\begin{align}\label{POINT1}
\sum_{|b|\leq 2}|Z^{b}u|\leq \nu_{0}.
\end{align}
Then we have
 \begin{align}
 |\partial_{\xi}A_1|&\leq C(|u_{\xi}|+|Zu_{\xi}|),~~|\partial_{\eta}A_1|\leq C(|u_{\eta}|+|Zu_{\eta}|),\\
 |\partial_{\xi}A_2|&\leq C(|u_{\eta}|+|Zu_{\eta}|),~~ |\partial_{\eta}A_2|\leq C(|u_{\eta}|+|Zu_{\eta}|),\\
 |\partial_{\xi}A_3|&\leq C(|u_{\xi}|+|Zu_{\xi}|),~~ |\partial_{\eta}A_3|\leq C(|u_{\xi}|+|Zu_{\xi}|),
 \end{align}
 where $C=C(\nu_0)$ is a positive constant depending on $\nu_0$.
 \end{Lemma}
 \begin{Lemma}\label{point222}
 Assume that $u$ is smooth, $A_1=A_1(u,u_{\xi},u_{\eta}), A_2=A_2(u,u_{\xi},u_{\eta})$, $A_3=A_3(u,u_{\xi},u_{\eta})$ and $F=F(u,u_{\xi},u_{\eta})$ satisfies \eqref{order1}, \eqref{order22},
\eqref{order33} and \eqref{order4}, respectively, and
\begin{align}\label{POINT12}
\sum_{|b|\leq 3}|Z^{b}u|\leq \nu_{1}.
\end{align}
Then for any multi-index $a, |a|=1, 2$, we have
\begin{align}
 |Z^{a}A_1|&\leq C\sum_{|b|\leq |a|}\big(|Z^{b}u_{\xi}|+|Z^{b}u_{\eta}|\big),\\
  |Z^{a}A_2|&\leq C\sum_{|b|\leq |a|}|Z^{b}u_{\eta}|,~
  |Z^{a}A_3|\leq C\sum_{|b|\leq |a|}|Z^{b}u_{\xi}|,\\
   |Z^{a}F|&\leq C\sum_{|b|+|c|\leq |a|}|Z^{b}u_{\xi}||Z^{c}u_{\eta}|.
 \end{align}
And for any multi-index $a,|a|= 3$, we also have
 \begin{align}
 |Z^{a}A_1|&\leq C\sum_{|b|\leq 3}\big(|Z^{b}u_{\xi}|+|Z^{b}u_{\eta}|\big),\\
  |Z^{a}A_2|&\leq C\sum_{|b|\leq 3}|Z^{b}u_{\eta}|+C|u_{\eta}|\sum_{|b|\leq 3}|Z^{b}u_{\xi}|,\\
  |Z^{a}A_3|&\leq C\sum_{|b|\leq 3}|Z^{b}u_{\xi}|+C|u_{\xi}|\sum_{|b|\leq 3}|Z^{b}u_{\eta}|,\\
   |Z^{a}F|&\leq C\sum_{|b|+|c|\leq 3}|Z^{b}u_{\xi}||Z^{c}u_{\eta}|.
 \end{align}
Here $C=C(\nu_1)$ is a positive constant depending on $\nu_1$.
 \end{Lemma}

 The following pointwise estimates will be used frequently in the sequel.
\begin{Lemma}\label{xuyao8DD8899}
Let $u$ be a smooth function on $\mathbb{R}^{+}\times \mathbb{R}^{+}$ with sufficient decay at spatial infinity, $k\geq 2$.
 Then it holds that
\begin{align}\label{224iijuio}
\|u(t,\cdot)\|_{L^{\infty}(\mathbb{R}^{+})}&\leq C
E^{1/2}(u(t)),
\end{align}
\begin{align}\label{xuyao217}
\sum_{|a|\leq k-2}\big(\|\langle\xi\rangle^{1+\delta}Z^{a}u_{\xi}\|_{L_{x}^{{\infty}}(\mathbb{R}^{+})}
+\|\langle \eta\rangle^{1+\delta}Z^{a}u_{\eta}\|_{L_{x}^{{\infty}}(\mathbb{R}^{+})}\big)&\leq C
E_k^{1/2}(u(t)),\\\label{xuddyao217}
\sum_{|a|\leq k-2}\big(\|\langle\xi\rangle^{3+3\delta}Z^{a}u_{\xi}\|_{L_{x}^{{\infty}}(\mathbb{R}^{+})}
+\|\langle \eta\rangle^{3+3\delta}Z^{a}u_{\eta}\|_{L_{x}^{{\infty}}(\mathbb{R}^{+})}\big)&\leq C
\mathbb{E}_k^{1/2}(u(t)),
\end{align}
and
\begin{align}\label{xuyao219}
\sum_{|a|\leq k-2}\big(\|\langle\eta\rangle^{-\frac{1+\delta}{2}}\langle\xi\rangle^{1+\delta}Z^{a}u_{\xi}\|_{L^2_{s}L_{x}^{{\infty}}}+\|\langle\xi\rangle^{-\frac{1+\delta}{2}}\langle \eta\rangle^{1+\delta}Z^{a}u_{\eta}\|_{L^2_{s}L_{x}^{{\infty}}}\big)
&\leq C
\mathcal {E}_k^{1/2}(u(t)),\\\label{xuyjuao219}
\sum_{|a|\leq k-2}\big(\|\langle\eta\rangle^{-\frac{1+\delta}{2}}\langle\xi\rangle^{3+3\delta}Z^{a}u_{\xi}\|_{L^2_{s}L_{x}^{{\infty}}}+\|\langle\xi\rangle^{-\frac{1+\delta}{2}}\langle \eta\rangle^{3+3\delta}Z^{a}u_{\eta}\|_{L^2_{s}L_{x}^{{\infty}}}\big)
&\leq C
\mathscr{E}_k^{1/2}(u(t)).
\end{align}
 \end{Lemma}
 \begin{proof}
  It follows from the fundamental theorem of calculus and H\"{o}lder inequality that
\begin{align}
&\|u(t,\cdot)\|_{L^{\infty}(\mathbb{R}^{+})}\leq \|u_x(t,\cdot)\|_{L_{x}^{1}(\mathbb{R}^{+})}\leq \|u_{\xi}\|_{L_{x}^{1}(\mathbb{R}^{+})}+\|u_{\eta}\|_{L_{x}^{1}(\mathbb{R}^{+})}\nonumber\\
&\leq \|\langle \xi\rangle^{-1-\delta}\|_{L_{x}^{2}(\mathbb{R}^{+})}\|\langle \xi\rangle^{1+\delta}u_{\xi}\|_{L_{x}^{2}(\mathbb{R}^{+})}+\|\langle \eta\rangle^{-1-\delta}\|_{L_{x}^{2}(\mathbb{R}^{+})}\|\langle \eta\rangle^{1+\delta}u_{\eta}\|_{L_{x}^{2}(\mathbb{R}^{+})}\nonumber\\
&\leq C\big(\|\langle \xi\rangle^{1+\delta}u_{\xi}\|_{L_{x}^{2}(\mathbb{R}^{+})}+\|\langle \eta\rangle^{1+\delta}u_{\eta}\|_{L_{x}^{2}(\mathbb{R}^{+})}\big)\leq C
E^{1/2}(u(t)).
\end{align}
For \eqref{xuyao217}, it can be proved by Sobolev embedding $H^1(\mathbb{R}^{+})\hookrightarrow L^{\infty}(\mathbb{R}^{+})$ and the following pointwise estimates
\begin{align}
|\partial_{x}\langle\xi\rangle^{1+\delta}|\leq C\langle\xi\rangle^{1+\delta},~~
|\partial_{x}\langle\eta\rangle^{1+\delta}|\leq C\langle\eta\rangle^{1+\delta}.
\end{align}
\eqref{xuddyao217} can be shown by the same way. While \eqref{xuyao219} is the consequence of Sobolev embedding $H^1(\mathbb{R}^{+})\hookrightarrow L^{\infty}(\mathbb{R}^{+})$ and the following fact
\begin{align}
|\partial_{x}(\langle\eta\rangle^{-\frac{1+\delta}{2}}\langle\xi\rangle^{1+\delta})|\leq C \langle\eta\rangle^{-\frac{1+\delta}{2}}\langle\xi\rangle^{1+\delta},~~
|\partial_{x}(\langle\xi\rangle^{-\frac{1+\delta}{2}}\langle\eta\rangle^{1+\delta})|\leq C \langle\xi\rangle^{-\frac{1+\delta}{2}}\langle\eta\rangle^{1+\delta}.
\end{align}
\eqref{xuyjuao219} can be proved similarly.
 \end{proof}

\subsection{Relationships of various energies}
The following two lemmas clarifies the relationships of various energies introduced in \S \ref{hjsedddffffffd}.
\begin{Lemma}\label{porpkoo}

We have
\begin{align}\label{bu1}
\overline{E}_k(u(t))&\leq C\big(\overline{\overline{E}}_k(u(t))+\widetilde{E}_k(u(t))\big),\\\label{bu1111}
\overline{\mathbb{E}}_k(u(t))&\leq C\big(\overline{\overline{\mathbb{E}}}_k(u(t))+\widetilde{\mathbb{E}}_k(u(t))\big),
\end{align}
and
\begin{align}\label{bu11}
\overline{\mathcal {E}}_k(u(t))&\leq C\big(\overline{\overline{\mathcal {E}}}_k(u(t))+\widetilde{\mathcal {E}}_k(u(t))\big),\\\label{bu1ddd1}
\overline{\mathscr{E}}_k(u(t))&\leq C\big(\overline{\overline{\mathscr{E}}}_k(u(t))+\widetilde{\mathscr{E}}_k(u(t))\big).
\end{align}

\end{Lemma}
\begin{proof}

In view of the definitions of energies, \eqref{bu1}, \eqref{bu1111}, \eqref{bu11} and \eqref{bu1ddd1} are consequences of the following pointwise estimates
\begin{align}\label{above1}
\sum_{l=0}^{k-1}|\partial_{\xi}^{l}u_{\xi}|\leq C\big( \sum_{|a|\leq k-2}|Z^{a}u_{\xi\eta}|+\sum_{l=0}^{k-1}|\partial_{t}^{l}u_{\xi}|\big),\\\label{above2}
\sum_{l=0}^{k-1}|\partial_{\eta}^{l}u_{\eta}|\leq C\big( \sum_{|a|\leq k-2}|Z^{a}u_{\xi\eta}|+\sum_{l=0}^{k-1}|\partial_{t}^{l}u_{\eta}|\big).
\end{align}
While~\eqref{above1} and \eqref{above2} can be shown by differential identity
\begin{align}
\partial_{\xi}= 2\partial_t -\partial_{\eta},~~(\text{i.e.,}~\partial_{\eta}= 2\partial_t -\partial_{\xi}  )
\end{align}
and inductive method.
\end{proof}

\begin{Lemma}\label{porpkoo1}
If $u$ satisfies \eqref{quasiwave}, $A_1, A_2$, $A_3$ and $F$ satisfies \eqref{order1}, \eqref{order22},
\eqref{order33} and \eqref{order4}, respectively, and
\begin{align}\label{tyuu78888}
\varepsilon_1=\sup_{0\leq s\leq t}E_4^{1/2}(u(s))
\end{align}
is sufficiently small.
Then we have
\begin{align}\label{bu2}
\sup_{0\leq s\leq t}\overline{\overline{E}}_4(u(s))&\leq C\sup_{0\leq s\leq t}{E}_4^2(u(s)),\\\label{bu21}
\sup_{0\leq s\leq t}\overline{\overline{\mathbb{E}}}_3(u(s))&\leq C\sup_{0\leq s\leq t}{E}_3(u(s))\sup_{0\leq s\leq t}{\mathbb{E}}_3(u(s)),
\end{align}
and
\begin{align}\label{bu22}
\overline{\overline{\mathcal {E}}}_4(u(t))&\leq C\sup_{0\leq s\leq t}{E}_4(u(s)){\mathcal {E}}_4(u(t)),\\\label{bu221}
\overline{\overline{\mathscr{E}}}_3(u(t))&\leq C\sup_{0\leq s\leq t}{E}_3(u(s)){\mathscr{E}}_3(u(t)).
\end{align}
\end{Lemma}
\begin{proof}
First, Lemma \ref{xuyao8DD8899} and \eqref{tyuu78888} imply \eqref{POINT12}, with $\nu_1 \simeq \varepsilon_1$, hence the conclusion of Lemma \ref{point222} holds throughout.
We only show \eqref{bu2} and \eqref{bu22}. While \eqref{bu21} and \eqref{bu221} can be proved similarly. In view of $u$ satisfies \eqref{quasiwave},
we get
\begin{align}
Z^{a}u_{\xi\eta}=\sum_{b+c=a}\lambda_{abc}\big( Z^{b}A_1Z^{c}u_{\xi\eta}+Z^{b}A_2Z^{c}u_{\xi\xi}+Z^{b}A_3Z^{c}u_{\eta\eta}\big)+Z^{a}F,
\end{align}
where $\lambda_{abc}$ are some constants. Thus for $|a|\leq 2$, it follows from Lemma \ref{point222} that
\begin{align}\label{xyaoddeeeeeedd}
|\langle\xi\rangle^{1+\delta}Z^{a}u_{\xi\eta}|&\leq C\sum_{|b|\leq 3}|Z^{b}u|\sum_{|c|\leq 3}|\langle\xi\rangle^{1+\delta}Z^{c}u_{\xi}|\nonumber\\
&+C\sum_{|b|\leq 2}|\langle\xi\rangle^{1+\delta}Z^{b}u_{\xi}|\sum_{|c|\leq 3}|Z^{c}u_{\eta}|,\\
|\langle\eta\rangle^{-\frac{1+\delta}{2}}\langle\xi\rangle^{1+\delta}Z^{a}u_{\xi\eta}|&\leq C\sum_{|b|\leq 3}|Z^{b}u|\sum_{|c|\leq 3}|\langle\eta\rangle^{-\frac{1+\delta}{2}}\langle\xi\rangle^{1+\delta}Z^{c}u_{\xi}|\nonumber\\
&+C\sum_{|b|\leq 2}|\langle\eta\rangle^{-\frac{1+\delta}{2}}\langle\xi\rangle^{1+\delta}Z^{b}u_{\xi}|\sum_{|c|\leq 3}|Z^{c}u_{\eta}|.
\end{align}
By \eqref{xyaoddeeeeeedd} and Lemma \ref{xuyao8DD8899} we have
\begin{align}\label{XUYASSSSO1}
&\|\langle\xi\rangle^{1+\delta}Z^{a}u_{\xi\eta}\|_{L_{x}^2(\mathbb{R}^{+})}\nonumber\\
&\leq C\sum_{|b|\leq 3}\|Z^{b}u\|_{L_{x}^{\infty}(\mathbb{R}^{+})}\sum_{|c|\leq 3}\|\langle\xi\rangle^{1+\delta}Z^{c}u_{\xi}\|_{L_{x}^2(\mathbb{R}^{+})}\nonumber\\
&+C\sum_{|b|\leq 2}\|\langle\xi\rangle^{1+\delta}Z^{b}u_{\xi}\|_{L_{x}^{\infty}(\mathbb{R}^{+})}\sum_{|c|\leq 3}\|Z^{c}u_{\eta}\|_{L_{x}^2(\mathbb{R}^{+})}\nonumber\\
&\leq C{E}_4(u(t))
\end{align}
and
\begin{align}\label{XUYASSDddddDDDDDDSSO1}
&\|\langle\eta\rangle^{-\frac{1+\delta}{2}}\langle\xi\rangle^{1+\delta}Z^{a}u_{\xi\eta}\|_{L_{s,x}^2}\nonumber\\
&\leq C\sum_{|b|\leq 3}\|Z^{b}u\|_{L_{s,x}^{\infty}}\sum_{|c|\leq 3}\|\langle\eta\rangle^{-\frac{1+\delta}{2}}\langle\xi\rangle^{1+\delta}Z^{c}u_{\xi}\|_{L_{s,x}^2}\nonumber\\
&+C\sum_{|b|\leq 2}\|\langle\eta\rangle^{-\frac{1+\delta}{2}}\langle\xi\rangle^{1+\delta}Z^{b}u_{\xi}\|_{L^{2}_{s}L_{x}^{\infty}}\sum_{|c|\leq 3}\|Z^{c}u_{\eta}\|_{L^{\infty}_{s}L_{x}^2}\nonumber\\
&\leq \sup_{0\leq s\leq t}{E}^{1/2}_4(u(s)){\mathcal {E}}^{1/2}_4(u(t)).
\end{align}
Similarly, we can also get
\begin{align}\label{XUYASSDDDSSO1}
&\|\langle\eta\rangle^{1+\delta}Z^{a}u_{\xi\eta}\|_{L_{x}^2(\mathbb{R}^{+})}\leq C{E}_4(u(t))
\end{align}
and
\begin{align}\label{XUYASSDDDDDDDSSO1}
&\|\langle\xi\rangle^{-\frac{1+\delta}{2}}\langle\eta\rangle^{1+\delta}Z^{a}u_{\xi\eta}\|_{L_{s,x}^2}\leq \sup_{0\leq s\leq t}{E}^{1/2}_4(u(s)){\mathcal {E}}^{1/2}_4(u(t)).
\end{align}
The combination of \eqref{XUYASSSSO1} and \eqref{XUYASSDDDSSO1}  implies \eqref{bu2}, and \eqref{bu22} is a consequence of \eqref{XUYASSDddddDDDDDDSSO1} and \eqref{XUYASSDDDDDDDSSO1}.
\end{proof}

\subsection{Energy estimates}
The following two lemmas on energy estimates are key points in the proof of Theorem \ref{mainthm}. The first one will be used in the control of high order energy.
\begin{Lemma}\label{xuoaj99hyyty}
Assume that $v: \mathbb{R}^{+}\times \mathbb{R}^{+}\longrightarrow \mathbb{R}^n$ satisfies the following system of linear wave equations
\begin{align}\label{timelike22666}
v_{\xi\eta}&=A_1(u,u_{\xi},u_{\eta})v_{\xi\eta}
+A_2(u,u_{\xi},u_{\eta})v_{\xi\xi}+A_3(u,u_{\xi},u_{\eta})v_{\eta\eta}+G,
\end{align}
and homogeneous Dirichlet boundary condition $v(t,0)=0$.
Here $A_1,A_2$, $A_3$ are symmetric and satisfies \eqref{order1}, \eqref{order22}, \eqref{order33}, respectively, $u$ and $G: \mathbb{R}^{+}\times \mathbb{R}^{+}\longrightarrow \mathbb{R}^n$ are some given vector valued functions of $(t,x)$, and
\begin{align}\label{tyuu78888}
\varepsilon_1=\sup_{0\leq s\leq t}\mathbb{E}^{1/2}_3(u(s))+\mathcal {E}^{1/2}_3(u(t))
\end{align}
is sufficiently small. Then it holds that\footnote{We denote the transpose of (column) vector $w\in \mathbb{R}^n$ by $w^{T}$.  }
\begin{align}\label{xiuppo999iu8}
&\sup_{0\leq s\leq t}{E}(v(s))+ {\mathcal{E}}(v(t))\nonumber\\
&\leq C E(v(0))+C\int_0^{t}\|\langle \xi\rangle^{2+2\delta}v^{T}_{\xi}G\|_{L_{x}^1(\mathbb{R}^{+})}ds+C\int_0^{t}\|\langle \eta\rangle^{2+2\delta}v^{T}_{\eta} G\|_{L_{x}^1(\mathbb{R}^{+})}ds.
\end{align}
\end{Lemma}
\begin{proof}
Multiply $2\psi(\eta)\phi(\xi)v^{{T}}_{\xi}$ on both sides of \eqref{timelike22666}. In view of the symmetry of $A_1, A_2$ and $A_3$, by Leibniz's rule we have
\begin{align}\label{lemkey1}
&\big(\psi(\eta)\phi(\xi)|v_{\xi}|^2\big)_{\eta}-\psi'(\eta)\phi(\xi)|v_{\xi}|^2\nonumber\\
&=\big(\psi(\eta)\phi(\xi)v_{\xi}^{{T}}A_1v_{\xi}\big)_{\eta}-\psi'(\eta)\phi(\xi)v_{\xi}^{{T}}A_1v_{\xi}
-\psi(\eta)\phi(\xi)v_{\xi}^{{T}}\partial_{\eta}A_1v_{\xi}\nonumber\\
&+\big(\psi(\eta)\phi(\xi)v_{\xi}^{{T}}A_2v_{\xi}\big)_{\xi}-\psi(\eta)\phi'(\xi)v_{\xi}^{{T}}A_2v_{\xi}
-\psi(\eta)\phi(\xi)v_{\xi}^{{T}}\partial_{\xi}A_2v_{\xi}\nonumber\\
&+\big(2\psi(\eta)\phi(\xi)v_{\xi}^{{T}}A_3v_{\eta}\big)_{\eta}-2\psi'(\eta)\phi(\xi)v_{\xi}^{{T}}A_3v_{\eta}
-2\psi(\eta)\phi(\xi)v_{\xi}^{{T}}\partial_{\eta}A_3v_{\eta}\nonumber\\
&-\big(\psi(\eta)\phi(\xi)v_{\eta}^{{T}}A_3v_{\eta}\big)_{\xi}+\psi(\eta)\phi'(\xi)v_{\eta}^{{T}}A_3v_{\eta}
+\psi(\eta)\phi(\xi)v_{\eta}^{{T}}\partial_{\xi}A_3v_{\eta}\nonumber\\
&+2\psi(\eta)\phi(\xi)v^{{T}}_{\xi}G.
\end{align}
Similarly, multiply $2\psi(\xi)\phi(\eta)v^{{T}}_{\eta}$ on both sides of \eqref{timelike22666}. The symmetry of $A_1, A_2$ and $A_3$ and Leibniz's rule also imply
\begin{align}\label{lemkey2}
&\big(\psi(\xi)\phi(\eta)|v_{\eta}|^2\big)_{\xi}-\psi'(\xi)\phi(\eta)|v_{\eta}|^2\nonumber\\
&=\big(\psi(\xi)\phi(\eta)v_{\eta}^{{T}}A_1v_{\eta}\big)_{\xi}-\psi'(\xi)\phi(\eta)v_{\eta}^{{T}}A_1v_{\eta}
-\psi(\xi)\phi(\eta)v_{\eta}^{{T}}\partial_{\xi}A_1v_{\eta}\nonumber\\
&+\big(2\psi(\xi)\phi(\eta)v_{\eta}^{{T}}A_2v_{\xi}\big)_{\xi}-2\psi'(\xi)\phi(\eta)v_{\eta}^{{T}}A_2v_{\xi}
-2\psi(\xi)\phi(\eta)v_{\eta}^{{T}}\partial_{\xi}A_2v_{\xi}\nonumber\\
&-\big(\psi(\xi)\phi(\eta)v_{\xi}^{{T}}A_2v_{\xi}\big)_{\eta}+\psi(\xi)\phi'(\eta)v_{\xi}^{{T}}A_2v_{\xi}
+\psi(\xi)\phi(\eta)v_{\xi}^{{T}}\partial_{\eta}A_2v_{\xi}\nonumber\\
&+\big(\psi(\xi)\phi(\eta)v_{\eta}^{{T}}A_3v_{\eta}\big)_{\eta}-\psi(\xi)\phi'(\eta)v_{\eta}^{{T}}A_3v_{\eta}
-\psi(\xi)\phi(\eta)v_{\eta}^{{T}}\partial_{\eta}A_3v_{\eta}\nonumber\\
&+2\psi(\xi)\phi(\eta)v^{{T}}_{\eta}G.
\end{align}

Integrating on $S_t=[0,t]\times\mathbb{R}^{+}$ on both sides of \eqref{lemkey1} and \eqref{lemkey2}, by the fundamental theorem of calculus we can get
\begin{align}\label{hj78999}
&\int_{\mathbb{R}^{+}}\big(e(t,x)+\widetilde{e}(t,x)\big)dx+\int_0^{t}\big(p(s,0)+\widetilde{p}(s,0)\big)ds+\int_0^{t}\!\!\int_{\mathbb{R}^{+}}q(s,x)dxds\nonumber\\
&=\int_{\mathbb{R}^{+}}\big(e(0,x)+\widetilde{e}(0,x)\big)dx+\int_0^{t}\!\!\int_{\mathbb{R}^{+}}\widetilde{q}(s,x)dxds\nonumber\\
&+2\int_0^{t}\!\!\int_{\mathbb{R}^{+}}\psi(\eta)\phi(\xi)v^{{T}}_{\xi}Gdxds+2\int_0^{t}\!\!\int_{\mathbb{R}^{+}}\psi(\xi)\phi(\eta)v^{{T}}_{\eta}Gdxds,
\end{align}
where
\begin{align}
e=\psi(\eta)\phi(\xi)|v_{\xi}|^2+\psi(\xi)\phi(\eta)|v_{\eta}|^2,
\end{align}
\begin{align}
\widetilde{e}&=-\psi(\eta)\phi(\xi)v_{\xi}^{{T}}A_1v_{\xi}-\psi(\eta)\phi(\xi)v_{\xi}^{{T}}A_2v_{\xi}
-2\psi(\eta)\phi(\xi)v_{\xi}^{{T}}A_3v_{\eta}+{\textcolor{blue}{\psi(\eta)\phi(\xi)v_{\eta}^{{T}}A_3v_{\eta}}}\nonumber\\
&-\psi(\xi)\phi(\eta)v_{\eta}^{{T}}A_1v_{\eta}
-2\psi(\xi)\phi(\eta)v_{\eta}^{{T}}A_2v_{\xi}
+{\textcolor{blue}{\psi(\xi)\phi(\eta)v_{\xi}^{{T}}A_2v_{\xi}}}
-\psi(\xi)\phi(\eta)v_{\eta}^{{T}}A_3v_{\eta},
\end{align}
\begin{align}\label{bianjie1}
p=\psi(\eta)\phi(\xi)|v_{\xi}|^2-\psi(\xi)\phi(\eta)|v_{\eta}|^2,
\end{align}
\begin{align}\label{bianjie2}
\widetilde{p}&=-\psi(\eta)\phi(\xi)v_{\xi}^{{T}}A_1v_{\xi}+\psi(\eta)\phi(\xi)v_{\xi}^{{T}}A_2v_{\xi}
-2\psi(\eta)\phi(\xi)v_{\xi}^{{T}}A_3v_{\eta}-\psi(\eta)\phi(\xi)v_{\eta}^{{T}}A_3v_{\eta}\nonumber\\
&+\psi(\xi)\phi(\eta)v_{\eta}^{{T}}A_1v_{\eta}
+2\psi(\xi)\phi(\eta)v_{\eta}^{{T}}A_2v_{\xi}
+\psi(\xi)\phi(\eta)v_{\xi}^{{T}}A_2v_{\xi}
-\psi(\xi)\phi(\eta)v_{\eta}^{{T}}A_3v_{\eta},
\end{align}
\begin{align}
q=-\psi'(\eta)\phi(\xi)|v_{\xi}|^2-\psi'(\xi)\phi(\eta)|v_{\eta}|^2,
\end{align}
and
\begin{align}\label{byw}
\widetilde{q}&=-\psi'(\eta)\phi(\xi)v_{\xi}^{{T}}A_1v_{\xi}-2\psi'(\eta)\phi(\xi)v_{\xi}^{{T}}A_3v_{\eta}
-\psi(\eta)\phi'(\xi)v_{\xi}^{{T}}A_2v_{\xi}+{\textcolor{blue}{\psi(\eta)\phi'(\xi)v_{\eta}^{{T}}A_3v_{\eta}}}
\nonumber\\
&-\psi(\eta)\phi(\xi)v_{\xi}^{{T}}\partial_{\eta}A_1v_{\xi}-\psi(\eta)\phi(\xi)v_{\xi}^{{T}}\partial_{\xi}A_2v_{\xi}
-2\psi(\eta)\phi(\xi)v_{\xi}^{{T}}\partial_{\eta}A_3v_{\eta}
+{\textcolor{blue}{\psi(\eta)\phi(\xi)v_{\eta}^{{T}}\partial_{\xi}A_3v_{\eta}}}
\nonumber\\
&-\psi'(\xi)\phi(\eta)v_{\eta}^{{T}}A_1v_{\eta}-2\psi'(\xi)\phi(\eta)v_{\eta}^{{T}}A_2v_{\xi}+{\textcolor{blue}{\psi(\xi)\phi'(\eta)v_{\xi}^{{T}}A_2v_{\xi}}}
-\psi(\xi)\phi'(\eta)v_{\eta}^{{T}}A_3v_{\eta}\nonumber\\
&-\psi(\xi)\phi(\eta)v_{\eta}^{{T}}\partial_{\xi}A_1v_{\eta}
-2\psi(\xi)\phi(\eta)v_{\eta}^{{T}}\partial_{\xi}A_2v_{\xi}
+{\textcolor{blue}{\psi(\xi)\phi(\eta)v_{\xi}^{{T}}\partial_{\eta}A_2v_{\xi}}}
-\psi(\xi)\phi(\eta)v_{\eta}^{{T}}\partial_{\eta}A_3v_{\eta}.
\end{align}


In view of \eqref{xuyao00099} and \eqref{rg56}, we have
\begin{align}
c^{-1}e(t,x)\leq |\langle \xi\rangle^{1+\delta}v_{\xi}|^2+|\langle \eta\rangle^{1+\delta}v_{\eta}|^2\leq ce(t,x)
\end{align}
It follows from \eqref{xuyao00099}, \eqref{rg56}, \eqref{order1}, \eqref{order22}, \eqref{order33} and Lemma \ref{xuyao8DD8899} that
\begin{align}
|\widetilde{e}(t,x)|&\leq C\langle \xi\rangle^{2+2\delta}|v_{\xi}^{{T}}A_1v_{\xi}|+C\langle \xi\rangle^{2+2\delta}|v_{\xi}^{{T}}A_2v_{\xi}|
+C\langle \xi\rangle^{2+2\delta}|v_{\xi}^{{T}}A_3v_{\eta}|+C\langle \xi\rangle^{2+2\delta}|v_{\eta}^{{T}}A_3v_{\eta}|\nonumber\\
&+C\langle \eta\rangle^{2+2\delta}|v_{\eta}^{{T}}A_1v_{\eta}|
+C\langle \eta\rangle^{2+2\delta}|v_{\eta}^{{T}}A_2v_{\xi}|
+C\langle \eta\rangle^{2+2\delta}|v_{\xi}^{{T}}A_2v_{\xi}|
+C\langle \eta\rangle^{2+2\delta}|v_{\eta}^{{T}}A_3v_{\eta}|\nonumber\\
&\leq C|\langle \xi\rangle^{1+\delta}v_{\xi}|^2(|u|+|u_{\xi}|+|u_{\eta}|)+C|\langle \xi\rangle^{1+\delta}v_{\xi}||\langle \xi\rangle^{1+\delta}u_{\xi}||v_{\eta}|+{\textcolor{blue}{C|\langle \xi\rangle^{2+2\delta}u_{\xi}||v_{\eta}|^2}}                                                       \nonumber\\
&+C|\langle \eta\rangle^{1+\delta}v_{\eta}|^2(|u|+|u_{\xi}|+|u_{\eta}|)
+C|\langle \eta\rangle^{1+\delta}v_{\eta}||\langle \eta\rangle^{1+\delta}u_{\eta}||v_{\xi}|+{\textcolor{blue}{C\langle \eta\rangle^{2+2\delta}u_{\eta}||v_{\xi}|^2}}
\nonumber\\
&\leq C\mathbb{E}_2^{1/2}(u(t))\big) \big(|\langle \xi\rangle^{1+\delta}v_{\xi}|^2+ |\langle \eta\rangle^{1+\delta}v_{\eta}|^2\big).
\end{align}
Noting \eqref{tyuu78888}, if $\varepsilon_1$ is sufficiently small, we can get
\begin{align}\label{siof00766hm}
\frac{c^{-1}}{2}\big(e(t,x)+\widetilde{e}(t,x)\big)\leq |\langle \xi\rangle^{1+\delta}v_{\xi}|^2+|\langle \eta\rangle^{1+\delta}v_{\eta}|^2\leq \frac{c}{2}\big(e(t,x)+\widetilde{e}(t,x)\big).
\end{align}
Note that $v(t,0)=0$ implies
\begin{align}\label{xuioowddd}
v_{\xi}(t,0)=v_{t}(t,0)+v_{x}(t,0)=v_{x}(t,0),~~v_{\eta}(t,0)=v_{t}(t,0)-v_{x}(t,0)=-v_{x}(t,0).
\end{align}
In view of \eqref{bianjie1}, by \eqref{xuioowddd} we obtain
\begin{align}\label{ling1}
p(t,0)=\psi(\frac{t}{2})\phi(\frac{t}{2})\big(|v_{\xi}(t,0)|^2-|v_{\eta}(t,0)|^2\big)=\psi(\frac{t}{2})\phi(\frac{t}{2})\big(|v_{x}(t,0)|^2-|v_{x}(t,0)|^2\big)=0.
\end{align}
In view of \eqref{bianjie2},
by \eqref{xuioowddd} we also have
\begin{align}\label{ling2}
\widetilde{p}(t,0)&=\psi(\frac{t}{2})\phi(\frac{t}{2})\big(-v^{{T}}_{\xi}(t,0)A_1v_{\xi}(t,0)+v^{{T}}_{\eta}(t,0)A_1v_{\eta}(t,0)\big)\nonumber\\
&+\psi(\frac{t}{2})\phi(\frac{t}{2})\big(v^{{T}}_{\xi}(t,0)A_2v_{\xi}(t,0)+2v^{{T}}_{\eta}(t,0)A_2v_{\xi}(t,0)+v^{{T}}_{\xi}(t,0)A_2v_{\xi}(t,0)\big)\nonumber\\
&+\psi(\frac{t}{2})\phi(\frac{t}{2})\big(-2v^{{T}}_{\xi}(t,0)A_3v_{\eta}(t,0)-v^{{T}}_{\eta}(t,0)A_3v_{\eta}(t,0)-v^{{T}}_{\eta}(t,0)A_3v_{\eta}(t,0)\big)\nonumber\\
&=\psi(\frac{t}{2})\phi(\frac{t}{2})\big(-v^{{T}}_{x}(t,0)A_1v_{x}(t,0)+v^{{T}}_{x}(t,0)A_1v_{x}(t,0)\big)\nonumber\\
&+\psi(\frac{t}{2})\phi(\frac{t}{2})\big(v^{{T}}_{x}(t,0)A_2v_{x}(t,0)-2v^{{T}}_{x}(t,0)A_2v_{x}(t,0)+v^{{T}}_{x}(t,0)A_2v_{x}(t,0)\big)\nonumber\\
&+\psi(\frac{t}{2})\phi(\frac{t}{2})\big(2v^{{T}}_{x}(t,0)A_3v_{x}(t,0)-v^{{T}}_{x}(t,0)A_3v_{x}(t,0)-v^{{T}}_{x}(t,0)A_3v_{x}(t,0)\big)\nonumber\\
&=0.
\end{align}
By \eqref{xuyao00099} and \eqref{here456}, we can obtain
\begin{align}\label{opp9pppp}
c^{-1}q(t,x)\leq \langle\eta\rangle^{-(1+\delta)}\langle \xi\rangle^{2+2\delta}|v_{\xi}|^2+\langle\xi\rangle^{-(1+\delta)}\langle \eta\rangle^{2+2\delta}|v_{\eta}|^2\leq cq(t,x).
\end{align}
Now Lemma \ref{xuyao8DD8899} and \eqref{tyuu78888} imply \eqref{POINT1}, with $\nu_0 \simeq \varepsilon_1$, hence the conclusions of Lemma \ref{xuyaouioo999} hold throughout.
By \eqref{byw} and Lemma \ref{xuyaouioo999},  we have
\begin{align}\label{qguji}
|\widetilde{q}(t,x)|&\leq C\langle \eta\rangle^{-(1+\delta)}\langle \xi\rangle^{2+2\delta}|v_{\xi}|^2(|u|+|u_{\xi}|+|u_{\eta}|)+C\langle \xi\rangle^{2+2\delta}|v_{\xi}|^2(|u_{\eta}|+|Zu_{\eta}|)\nonumber\\
&+ C\langle \eta\rangle^{-(1+\delta)}\langle \xi\rangle^{2+2\delta}|v_{\xi}||v_{\eta}||u_{\xi}|+C\langle \xi\rangle^{2+2\delta}|v_{\xi}||v_{\eta}|(|u_{\xi}|+|Zu_{\xi}|)\nonumber\\
&+{\textcolor{blue}{C\langle \xi\rangle^{2+2\delta}|v_{\eta}|^2(|u_{\xi}|+|Zu_{\xi}|)}}\nonumber\\
&+ C\langle \xi\rangle^{-(1+\delta)}\langle \eta\rangle^{2+2\delta}|v_{\eta}|^2(|u|+|u_{\xi}|+|u_{\eta}|)+C\langle\eta\rangle^{2+2\delta}|v_{\eta}|^2(|u_{\xi}|+|Zu_{\xi}|)\nonumber\\
&+C\langle \xi\rangle^{-(1+\delta)}\langle \eta\rangle^{2+2\delta}|v_{\xi}||v_{\eta}||u_{\eta}|+C\langle \eta\rangle^{2+2\delta}|v_{\xi}||v_{\eta}|(|u_{\eta}|+|Zu_{\eta}|)\nonumber\\
&+C{\textcolor{blue}{\langle \eta\rangle^{2+2\delta}|v_{\xi}|^2(|u_{\eta}|+|Zu_{\eta}|)}}\nonumber\\
&\leq C|\langle \eta\rangle^{-\frac{1+\delta}{2}}\langle \xi\rangle^{1+\delta}v_{\xi}|^2(|u|+|u_{\xi}|+\langle \eta\rangle^{1+\delta}|u_{\eta}|+\langle \eta\rangle^{1+\delta}|Zu_{\eta}|)\nonumber\\
&+C|\langle \eta\rangle^{-\frac{1+\delta}{2}}\langle \xi\rangle^{1+\delta}v_{\xi}||\langle \eta\rangle^{-\frac{1+\delta}{2}}\langle \xi\rangle^{1+\delta}(|u_{\xi}|+|Zu_{\xi}|)|\langle \eta\rangle^{1+\delta}v_{\eta}|\nonumber\\
&+C{\textcolor{blue}{\langle \xi\rangle^{3+3\delta}(|u_{\xi}|+|Zu_{\xi}|)}}|\langle \xi\rangle^{-\frac{1+\delta}{2}}\langle \eta\rangle^{1+\delta}v_{\eta}|^2\nonumber\\
&+C|\langle \xi\rangle^{-\frac{1+\delta}{2}}\langle \eta\rangle^{1+\delta}v_{\eta}|^2(|u|+\langle \xi\rangle^{1+\delta}|u_{\xi}|+\langle \xi\rangle^{1+\delta}|Zu_{\xi}|+|u_{\eta}|)\nonumber\\
&+C|\langle \xi\rangle^{-\frac{1+\delta}{2}}\langle \eta\rangle^{1+\delta}v_{\eta}||\langle \xi\rangle^{-\frac{1+\delta}{2}}\langle \eta\rangle^{1+\delta}(|u_{\eta}|+|Zu_{\eta}|)|\langle \xi\rangle^{1+\delta}v_{\xi}|\nonumber\\
&+C{\textcolor{blue}{\langle \eta\rangle^{3+3\delta}(|u_{\eta}|+|Zu_{\eta}|)|}}|\langle \eta\rangle^{-\frac{1+\delta}{2}}\langle \xi\rangle^{1+\delta}v_{\xi}|^2.
\end{align}

The combination of \eqref{hj78999}, \eqref{siof00766hm}, \eqref{ling1}, \eqref{ling2} and \eqref{opp9pppp} implies
\begin{align}\label{xddd89uii}
&\sup_{0\leq s\leq t}E(v(s))+ \mathcal{E}(v(t))\nonumber\\
&\leq C E(v(0))+C\int_0^{t}\|\widetilde{q}(s,\cdot)\|_{L^1(\mathbb{R}^{+})}ds\nonumber\\
&+C\int_0^{t}\|\langle \xi\rangle^{2+2\delta}v^{T}_{\xi}G\|_{L_{x}^1(\mathbb{R}^{+})}ds+C\int_0^{t}\|\langle \eta\rangle^{2+2\delta}v^{T}_{\eta} G\|_{L_{x}^1(\mathbb{R}^{+})}ds.
\end{align}
It follows from \eqref{qguji}, H\"{o}lder inequality and Lemma \ref{xuyao8DD8899} that
\begin{align}\label{qgu999ji}
&\int_0^{t}\|\widetilde{q}(s,\cdot)\|_{L^1(\mathbb{R}^{+})}ds\nonumber\\
&\leq C\|\langle \eta\rangle^{-\frac{1+\delta}{2}}\langle \xi\rangle^{1+\delta}v_{\xi}\|_{L^2_{s,x}}^2(\|u\|_{L^{\infty}_{s,x}}+\|u_{\xi}\|_{L^{\infty}_{s,x}}+\|\langle \eta\rangle^{1+\delta}|u_{\eta}|+\langle \eta\rangle^{1+\delta}|Zu_{\eta}|\|_{L^{\infty}_{s,x}})\nonumber\\
&+C\|\langle \eta\rangle^{-\frac{1+\delta}{2}}\langle \xi\rangle^{1+\delta}v_{\xi}\|_{L^{2}_{s,x}}\|\langle \eta\rangle^{-\frac{1+\delta}{2}}\langle \xi\rangle^{1+\delta}(|u_{\xi}|+|Zu_{\xi}|)\|_{L^{2}_{s}L^{\infty}_{x}}\|\langle \eta\rangle^{1+\delta}v_{\eta}\|_{L^{\infty}_{s}L^{2}_{x}}\nonumber\\
&+C{\textcolor{blue}{\|\langle \xi\rangle^{3+3\delta}(|u_{\xi}|+|Zu_{\xi}|)\|_{L^{\infty}_{s,x}}}}\|\langle \xi\rangle^{-\frac{1+\delta}{2}}\langle \eta\rangle^{1+\delta}v_{\eta}\|_{L^2_{s,x}}^2\nonumber\\
&+C\|\langle \xi\rangle^{-\frac{1+\delta}{2}}\langle \eta\rangle^{1+\delta}v_{\eta}\|_{L^{2}_{s,x}}^2(\|u\|_{L^{\infty}_{s,x}}+\|\langle \xi\rangle^{1+\delta}|u_{\xi}|+\langle \xi\rangle^{1+\delta}|Zu_{\xi}|\|_{L^{\infty}_{s,x}}+\|u_{\eta}\|_{L^{\infty}_{s,x}})\nonumber\\
&+C\|\langle \xi\rangle^{-\frac{1+\delta}{2}}\langle \eta\rangle^{1+\delta}v_{\eta}\|_{L^{2}_{s,x}}\|\langle \xi\rangle^{-\frac{1+\delta}{2}}\langle \eta\rangle^{1+\delta}(|u_{\eta}|+|Zu_{\eta}|)\|_{L^{2}_{s}L^{\infty}_{x}}\|\langle \xi\rangle^{1+\delta}v_{\xi}\|_{L^{\infty}_{s}L^{2}_{x}}\nonumber\\
&+C{\textcolor{blue}{\|\langle \eta\rangle^{3+3\delta}(|u_{\eta}|+|Zu_{\eta}|)\|_{L^{\infty}_{s,x}}}}\|\langle \eta\rangle^{-\frac{1+\delta}{2}}\langle \xi\rangle^{1+\delta}v_{\xi}\|_{L^2_{s,x}}^2\nonumber\\
&\leq C{\textcolor{blue}{\sup_{0\leq s\leq t}\mathbb{E}_3^{1/2}(u(s))}}\mathcal {E}(v(t))+C\mathcal {E}_3^{1/2}(u(t))\mathcal {E}^{1/2}(v(t))\sup_{0\leq s\leq t}{E}^{1/2}(v(s)).
\end{align}

Finally, by \eqref{tyuu78888}, \eqref{xddd89uii} and \eqref{qgu999ji}, we have
\begin{align}\label{xddd89ui999i}
&\sup_{0\leq s\leq t}E(v(s))+ \mathcal{E}(v(t))\nonumber\\
&\leq C E(v(0))+C\big(\sup_{0\leq s\leq t}\mathbb{E}_3^{1/2}(u(s))+\mathcal {E}_3^{1/2}(u(t))\big)\big(\sup_{0\leq s\leq t}E(v(s))+ \mathcal{E}(v(t))\big)\nonumber\\
&+C\int_0^{t}\|\langle \xi\rangle^{2+2\delta}v^{T}_{\xi} G\|_{L_{x}^1(\mathbb{R}^{+})}ds+C\int_0^{t}\|\langle \eta\rangle^{2+2\delta}v^{T}_{\eta} G\|_{L_{x}^1(\mathbb{R}^{+})}ds\nonumber\\
&\leq C E(v(0))+C\varepsilon_1\big(\sup_{0\leq s\leq t}E(v(s))+ \mathcal{E}(v(t))\big)\nonumber\\
&+C\int_0^{t}\|\langle \xi\rangle^{2+2\delta}v^{T}_{\xi} G\|_{L_{x}^1(\mathbb{R}^{+})}ds+C\int_0^{t}\|\langle \eta\rangle^{2+2\delta}v^{T}_{\eta} G\|_{L_{x}^1(\mathbb{R}^{+})}ds.
\end{align}
If $\varepsilon_1$ is sufficiently small, we can get \eqref{xiuppo999iu8}.
\end{proof}

 The following lemma will be used in the control of low order energy in the proof of Theorem \ref{mainthm}.
\begin{Lemma}\label{xuoaj99hyyty9}
Assume that $v: \mathbb{R}^{+}\times \mathbb{R}^{+}\longrightarrow \mathbb{R}^n$ satisfies the following system of linear wave equations
\begin{align}\label{timelike226669}
v_{\xi\eta}&=A_1(u,u_{\xi},u_{\eta})v_{\xi\eta}
+A_2(u,u_{\xi},u_{\eta})v_{\xi\xi}+A_3(u,u_{\xi},u_{\eta})v_{\eta\eta}+G,
\end{align}
and homogeneous Dirichlet boundary condition $v(t,0)=0$.
Here $A_1,A_2$, $A_3$ are symmetric and satisfies \eqref{order1}, \eqref{order22}, \eqref{order33}, respectively, $u$ and $G: \mathbb{R}^{+}\times \mathbb{R}^{+}\longrightarrow \mathbb{R}^n$ are some given vector valued functions of $(t,x)$, and
\begin{align}\label{tyuu788889}
\varepsilon_1=\sup_{0\leq s\leq t}{E}^{1/2}_3(u(s))
\end{align}
is sufficiently small. Then it holds that
\begin{align}\label{xiuppo999iu89}
&\sup_{0\leq s\leq t}\mathbb{E}(v(s))+\mathscr{E}(v(t))\nonumber\\
&\leq C \mathbb{E}(v(0))+C\mathscr{E}_2^{1/2}(u(t))\mathscr{E}^{1/2}(v(t))\sup_{0\leq s\leq t}E^{1/2}_2(v(t))\nonumber\\
&+C\int_0^{t}\|\langle \xi\rangle^{6+6\delta}v^{T}_{\xi}G\|_{L_{x}^1(\mathbb{R}^{+})}ds+C\int_0^{t}\|\langle \eta\rangle^{6+6\delta}v^{T}_{\eta} G\|_{L_{x}^1(\mathbb{R}^{+})}ds.
\end{align}
\end{Lemma}
\begin{proof}
Multiply $2\psi(\eta)\theta(\xi)v^{{T}}_{\xi}$ on both sides of \eqref{timelike226669}. Noting the symmetry of $A_1, A_2$ and $A_3$, by Leibniz's rule we have
\begin{align}\label{lemkey19}
&\big(\psi(\eta)\theta(\xi)|v_{\xi}|^2\big)_{\eta}-\psi'(\eta)\theta(\xi)|v_{\xi}|^2\nonumber\\
&=\big(\psi(\eta)\theta(\xi)v_{\xi}^{{T}}A_1v_{\xi}\big)_{\eta}-\psi'(\eta)\theta(\xi)v_{\xi}^{{T}}A_1v_{\xi}
-\psi(\eta)\theta(\xi)v_{\xi}^{{T}}\partial_{\eta}A_1v_{\xi}\nonumber\\
&+\big(\psi(\eta)\theta(\xi)v_{\xi}^{{T}}A_2v_{\xi}\big)_{\xi}-\psi(\eta)\theta'(\xi)v_{\xi}^{{T}}A_2v_{\xi}
-\psi(\eta)\theta(\xi)v_{\xi}^{{T}}\partial_{\xi}A_2v_{\xi}\nonumber\\
&+2\psi(\eta)\theta(\xi)v^{{T}}_{\xi}A_3v_{\eta\eta}+2\psi(\eta)\theta(\xi)v^{{T}}_{\xi}G.
\end{align}
Now, multiply $\psi(\xi)\theta(\eta)v^{{T}}_{\eta}$ on both sides of \eqref{timelike226669}. The symmetry of $A_1, A_2$ and $A_3$ and Leibniz's rule also imply
\begin{align}\label{lemkey29}
&\big({\textcolor{blue}{\frac{1}{2}}}\psi(\xi)\theta(\eta)|v_{\eta}|^2\big)_{\xi}-\frac{1}{2}\psi'(\xi)\theta(\eta)|v_{\eta}|^2\nonumber\\
&=\big(\frac{1}{2}\psi(\xi)\theta(\eta)v_{\eta}^{{T}}A_1v_{\eta}\big)_{\xi}-\frac{1}{2}\psi'(\xi)\theta(\eta)v_{\eta}^{{T}}A_1v_{\eta}
-\frac{1}{2}\psi(\xi)\theta(\eta)v_{\eta}^{{T}}\partial_{\xi}A_1v_{\eta}
+\psi(\xi)\theta(\eta)v^{{T}}_{\eta}A_2v_{\xi\xi}\nonumber\\
&+\frac{1}{2}\big(\psi(\xi)\theta(\eta)v^{{T}}_{\eta}A_3v_{\eta}\big)_{\eta}-\frac{1}{2}\psi(\xi)\theta'(\eta)v_{\eta}^{{T}}A_3v_{\eta}
-\frac{1}{2}\psi(\xi)\theta(\eta)v^{{T}}_{\eta}\partial_{\eta}A_3v_{\eta}\nonumber\\
&+\psi(\xi)\theta(\eta)v^{{T}}_{\eta}G.
\end{align}

Integrating on $S_t=[0,t]\times\mathbb{R}^{+}$ on both sides of \eqref{lemkey19} and \eqref{lemkey29}, by the fundamental theorem of calculus we have
\begin{align}\label{hj789999}
&\int_{\mathbb{R}^{+}}\big(e(t,x)+\widetilde{e}(t,x)\big)dx+\int_0^{t}\big(p(s,0)+\widetilde{p}(s,0)\big)ds+\int_0^{t}\!\!\int_{\mathbb{R}^{+}}q(s,x)dxds
\nonumber\\
&=\int_{\mathbb{R}^{+}}\big(e(0,x)+\widetilde{e}(0,x)\big)dx+\int_0^{t}\!\!\int_{\mathbb{R}^{+}}\widetilde{q}(s,x)dxds\nonumber\\
&+2\int_0^{t}\!\!\int_{\mathbb{R}^{+}}\psi(\eta)\theta(\xi)v^{{T}}_{\xi}A_3v_{\eta\eta}dxds+
\int_0^{t}\!\!\int_{\mathbb{R}^{+}}\psi(\xi)\theta(\eta)v^{{T}}_{\eta}A_2v_{\xi\xi}dxds\nonumber\\
&+2\int_0^{t}\!\!\int_{\mathbb{R}^{+}}\psi(\eta)\theta(\xi)v^{{T}}_{\xi}Gdxds+\int_0^{t}\!\!\int_{\mathbb{R}^{+}}\psi(\xi)\theta(\eta)v^{{T}}_{\eta}Gdxds,
\end{align}
where
\begin{align}
e=\psi(\eta)\theta(\xi)|v_{\xi}|^2+\frac{1}{2}\psi(\xi)\theta(\eta)|v_{\eta}|^2,
\end{align}
\begin{align}
\widetilde{e}&=-\psi(\eta)\theta(\xi)v_{\xi}^{{T}}A_1v_{\xi}-\psi(\eta)\theta(\xi)v_{\xi}^{{T}}A_2v_{\xi}
\nonumber\\
&-\frac{1}{2}\psi(\xi)\theta(\eta)v_{\eta}^{{T}}A_1v_{\eta}
-\frac{1}{2}\psi(\xi)\theta(\eta)v_{\eta}^{{T}}A_3v_{\eta},
\end{align}
\begin{align}\label{bianjie19}
p=\psi(\eta)\theta(\xi)|v_{\xi}|^2-{\textcolor{blue}{\frac{1}{2}}}\psi(\xi)\theta(\eta)|v_{\eta}|^2,
\end{align}
\begin{align}\label{bianjie29}
\widetilde{p}&=-\psi(\eta)\theta(\xi)v_{\xi}^{{T}}A_1v_{\xi}+\psi(\eta)\theta(\xi)v_{\xi}^{{T}}A_2v_{\xi}
\nonumber\\
&+\frac{1}{2}\psi(\xi)\theta(\eta)v_{\eta}^{{T}}A_1v_{\eta}
-\frac{1}{2}\psi(\xi)\theta(\eta)v_{\eta}^{{T}}A_3v_{\eta},
\end{align}
\begin{align}
q=-\psi'(\eta)\theta(\xi)|v_{\xi}|^2-\frac{1}{2}\psi'(\xi)\theta(\eta)|v_{\eta}|^2,
\end{align}
and
\begin{align}\label{byw9}
\widetilde{q}&=-\psi'(\eta)\theta(\xi)v_{\xi}^{{T}}A_1v_{\xi}
-\psi(\eta)\theta'(\xi)v_{\xi}^{{T}}A_2v_{\xi}
\nonumber\\
&-\psi(\eta)\theta(\xi)v_{\xi}^{{T}}\partial_{\eta}A_1v_{\xi}-\psi(\eta)\theta(\xi)v_{\xi}^{{T}}\partial_{\xi}A_2v_{\xi}
\nonumber\\
&-\frac{1}{2}\psi'(\xi)\theta(\eta)v_{\eta}^{{T}}A_1v_{\eta}
-\frac{1}{2}\psi(\xi)\theta'(\eta)v_{\eta}^{{T}}A_3v_{\eta}\nonumber\\
&-\frac{1}{2}\psi(\xi)\theta(\eta)v_{\eta}^{{T}}\partial_{\xi}A_1v_{\eta}
-\frac{1}{2}\psi(\xi)\theta(\eta)v_{\eta}^{{T}}\partial_{\eta}A_3v_{\eta}.
\end{align}


In view of \eqref{xuyao00099} and \eqref{rg56}, we have
\begin{align}
c^{-1}e(t,x)\leq |\langle \xi\rangle^{3+3\delta}v_{\xi}|^2+\frac{1}{2}|\langle \eta\rangle^{3+3\delta}v_{\eta}|^2\leq ce(t,x)
\end{align}
It follows from \eqref{xuyao00099}, \eqref{rg56}, \eqref{order1}, \eqref{order22}, \eqref{order33} and Lemma \ref{xuyao8DD8899} that
\begin{align}
|\widetilde{e}(t,x)|&\leq C\langle \xi\rangle^{6+6\delta}|v_{\xi}^{{T}}A_1v_{\xi}|+C\langle \xi\rangle^{6+6\delta}|v_{\xi}^{{T}}A_2v_{\xi}|
\nonumber\\
&+C\langle \eta\rangle^{6+6\delta}|v_{\eta}^{{T}}A_1v_{\eta}|
+C\langle \eta\rangle^{6+6\delta}|v_{\eta}^{{T}}A_3v_{\eta}|\nonumber\\
&\leq C \big(|\langle \xi\rangle^{3+3\delta}v_{\xi}|^2+ |\langle \eta\rangle^{3+3\delta}v_{\eta}|^2\big)(|u|+|u_{\xi}|+|u_{\eta}|)
\nonumber\\
&\leq C{E}_2^{1/2}(u(t)) \big(|\langle \xi\rangle^{3+3\delta}v_{\xi}|^2+ |\langle \eta\rangle^{3+3\delta}v_{\eta}|^2\big).
\end{align}
Noting \eqref{tyuu788889}, if $\varepsilon_1$ is sufficiently small, we can get
\begin{align}\label{siof00766hm9}
\frac{c^{-1}}{2}\big(e(t,x)+\widetilde{e}(t,x)\big)\leq |\langle \xi\rangle^{3+3\delta}v_{\xi}|^2+\frac{1}{2}|\langle \eta\rangle^{3+3\delta}v_{\eta}|^2\leq \frac{c}{2}\big(e(t,x)+\widetilde{e}(t,x)\big).
\end{align}
Note that $v(t,0)=0$ implies
\begin{align}\label{xuioowddd9}
v_{\xi}(t,0)=v_{t}(t,0)+v_{x}(t,0)=v_{x}(t,0),~~v_{\eta}(t,0)=v_{t}(t,0)-v_{x}(t,0)=-v_{x}(t,0).
\end{align}
In view of \eqref{bianjie19}, by \eqref{xuioowddd9} we have
\begin{align}\label{ling19}
p(t,0)=\psi(\frac{t}{2})\theta(\frac{t}{2})\big(|v_{\xi}(t,0)|^2-\frac{1}{2}|v_{\eta}(t,0)|^2\big)=\frac{1}{2}\psi(\frac{t}{2})\theta(\frac{t}{2})|v_{x}(t,0)|^2.
\end{align}
In view of \eqref{bianjie29},
by \eqref{xuioowddd9} we also have
\begin{align}\label{ling29}
\widetilde{p}(t,0)&=\psi(\frac{t}{2})\theta(\frac{t}{2})\big(-v^{{T}}_{\xi}(t,0)A_1v_{\xi}(t,0)+\frac{1}{2}v^{{T}}_{\eta}(t,0)A_1v_{\eta}(t,0)\big)\nonumber\\
&+\psi(\frac{t}{2})\theta(\frac{t}{2})\big(v^{{T}}_{\xi}(t,0)A_2v_{\xi}(t,0)-\frac{1}{2}v^{{T}}_{\eta}(t,0)A_3v_{\eta}(t,0)\big)\nonumber\\
&=\psi(\frac{t}{2})\theta(\frac{t}{2})v^{{T}}_{x}(t,0)(-\frac{1}{2}A_1+A_2-\frac{1}{2}A_3)v_{x}(t,0).
\end{align}
Thus by \eqref{order1}, \eqref{order22}, \eqref{order33} and Lemma \ref{xuyao8DD8899} we can get
\begin{align}\label{ling992}
|\widetilde{p}(t,0)|\leq C\psi(\frac{t}{2})\theta(\frac{t}{2}){E}_2^{1/2}(u(t))|v_{x}(t,0)|^2.
\end{align}
It follows from \eqref{ling19} and \eqref{ling992} that if $\varepsilon_1$ is sufficiently small, we can get
\begin{align}\label{ling9967772}
{p}(t,0)+\widetilde{p}(t,0)\geq \frac{1}{4}\psi(\frac{t}{2})\theta(\frac{t}{2})|v_{x}(t,0)|^2\geq 0.
\end{align}
By \eqref{xuyao00099} and \eqref{here456}, we can obtain
\begin{align}\label{opp9pppp9}
c^{-1}q(t,x)\leq \langle\eta\rangle^{-(1+\delta)}\langle \xi\rangle^{6+6\delta}|v_{\xi}|^2+\langle\xi\rangle^{-(1+\delta)}\langle \eta\rangle^{6+6\delta}|v_{\eta}|^2\leq cq(t,x).
\end{align}
Now Lemma \ref{xuyao8DD8899} and \eqref{tyuu788889} imply \eqref{POINT1}, with $\nu_0 \simeq \varepsilon_1$, hence the conclusions of Lemma \ref{xuyaouioo999} hold throughout.
By \eqref{byw9} and Lemma \ref{xuyaouioo999},  we have
\begin{align}\label{qguji9}
|\widetilde{q}(t,x)|&\leq C\langle \eta\rangle^{-(1+\delta)}\langle \xi\rangle^{6+6\delta}|v_{\xi}|^2(|u|+|u_{\xi}|+|u_{\eta}|)+C\langle \xi\rangle^{6+6\delta}|v_{\xi}|^2(|u_{\eta}|+|Zu_{\eta}|)\nonumber\\
&+ C\langle \xi\rangle^{-(1+\delta)}\langle \eta\rangle^{6+6\delta}|v_{\eta}|^2(|u|+|u_{\xi}|+|u_{\eta}|)+C\langle\eta\rangle^{6+6\delta}|v_{\eta}|^2(|u_{\xi}|+|Zu_{\xi}|)\nonumber\\
&\leq C|\langle \eta\rangle^{-\frac{1+\delta}{2}}\langle \xi\rangle^{3+3\delta}v_{\xi}|^2(|u|+|u_{\xi}|+\langle \eta\rangle^{1+\delta}|u_{\eta}|+\langle \eta\rangle^{1+\delta}|Zu_{\eta}|)\nonumber\\
&+C|\langle \xi\rangle^{-\frac{1+\delta}{2}}\langle \eta\rangle^{3+3\delta}v_{\eta}|^2(|u|+\langle \xi\rangle^{1+\delta}|u_{\xi}|+\langle \xi\rangle^{1+\delta}|Zu_{\xi}|+|u_{\eta}|).
\end{align}

The combination of \eqref{hj789999}, \eqref{siof00766hm9}, \eqref{ling9967772} and \eqref{opp9pppp9} implies
\begin{align}\label{xddd89uii9}
&\sup_{0\leq s\leq t}\mathbb{E}(v(s))+\mathscr{E}(v(t))\nonumber\\
&\leq C \mathbb{E}(v(0))+C\int_0^{t}\|\widetilde{q}(s,\cdot)\|_{L^1(\mathbb{R}^{+})}ds\nonumber\\
&+C\int_0^{t}\big(\|\langle \xi\rangle^{6+6\delta}u_{\xi}v_{\xi}v_{\eta\eta}\|_{L^1(\mathbb{R}^{+})} +\|\langle \eta\rangle^{6+6\delta}u_{\eta}v_{\eta}v_{\xi\xi}\|_{L^1(\mathbb{R}^{+})}  \big)ds\nonumber\\
&+C\int_0^{t}\|\langle \xi\rangle^{6+6\delta}v^{T}_{\xi}G\|_{L_{x}^1(\mathbb{R}^{+})}ds+C\int_0^{t}\|\langle \eta\rangle^{6+6\delta}v^{T}_{\eta} G\|_{L_{x}^1(\mathbb{R}^{+})}ds.
\end{align}
It follows from \eqref{qguji9}, H\"{o}lder inequality and Lemma \ref{xuyao8DD8899} that
\begin{align}\label{qgu999ji9}
&\int_0^{t}\|\widetilde{q}(s,\cdot)\|_{L^1(\mathbb{R}^{+})}ds\nonumber\\
&\leq C\|\langle \eta\rangle^{-\frac{1+\delta}{2}}\langle \xi\rangle^{3+3\delta}v_{\xi}\|_{L^2_{s,x}}^2(\|u\|_{L^{\infty}_{s,x}}+\|u_{\xi}\|_{L^{\infty}_{s,x}}+\|\langle \eta\rangle^{1+\delta}(|u_{\eta}|+|Zu_{\eta}|)\|_{L^{\infty}_{s,x}})\nonumber\\
&+ C\|\langle \xi\rangle^{-\frac{1+\delta}{2}}\langle \eta\rangle^{3+3\delta}v_{\eta}\|_{L^2_{s,x}}^2(\|u\|_{L^{\infty}_{s,x}}+\|\langle \xi\rangle^{1+\delta}(|u_{\xi}|+|Zu_{\xi}|)\|_{L^{\infty}_{s,x}}+\|u_{\eta}\|_{L^{\infty}_{s,x}})\nonumber\\
&\leq C\sup_{0\leq s\leq t}E_3^{1/2}(u(s))\mathscr{E}(v(t)).
\end{align}
By H\"{o}lder inequality and Lemma \ref{xuyao8DD8899}, we also have
\begin{align}\label{xddd89uiiddddd9}
&\int_0^{t}\big(\|\langle \xi\rangle^{6+6\delta}u_{\xi}v_{\xi}v_{\eta\eta}\|_{L^1(\mathbb{R}^{+})} +\|\langle \eta\rangle^{6+6\delta}u_{\eta}v_{\eta}v_{\xi\xi}\|_{L^1(\mathbb{R}^{+})}  \big)ds\nonumber\\
&\leq \|\langle \eta\rangle^{-\frac{1+\delta}{2}}\langle \xi\rangle^{3+3\delta}v_{\xi}\|_{L^2_{s,x}}\|\langle \eta\rangle^{-\frac{1+\delta}{2}}\langle \xi\rangle^{3+3\delta}u_{\xi}\|_{L^{2}_{s}L^{\infty}_{s}}\|\langle \eta\rangle^{1+\delta}v_{\eta\eta}\|_{L^{\infty}_{s}L^2_{x}}\nonumber\\
&+\|\langle \xi\rangle^{-\frac{1+\delta}{2}}\langle \eta\rangle^{3+3\delta}v_{\eta}\|_{L^2_{s,x}}\|\langle \xi\rangle^{-\frac{1+\delta}{2}}\langle \eta\rangle^{3+3\delta}u_{\eta}\|_{L^{2}_{s}L^{\infty}_{s}}\|\langle \xi\rangle^{1+\delta}v_{\xi\xi}\|_{L^{\infty}_{s}L^2_{x}}\nonumber\\
&\leq C\mathscr{E}_2^{1/2}(u(t))\mathscr{E}^{1/2}(v(t))\sup_{0\leq s\leq t}E^{1/2}_2(v(t)).
\end{align}

Finally, by \eqref{xddd89uii9}, \eqref{qgu999ji9} and \eqref{xddd89uiiddddd9}, we have
\begin{align}\label{xddd89ui999i9}
&\sup_{0\leq s\leq t}\mathbb{E}(v(s))+\mathscr{E}(v(t))\nonumber\\
&\leq C \mathbb{E}(v(0))+C\sup_{0\leq s\leq t}E_3^{1/2}(u(s))\mathscr{E}(v(t))+C\mathscr{E}_2^{1/2}(u(t))\mathscr{E}^{1/2}(v(t))\sup_{0\leq s\leq t}E^{1/2}_2(v(t))\nonumber\\
&+C\int_0^{t}\|\langle \xi\rangle^{6+6\delta}v^{T}_{\xi}G\|_{L_{x}^1(\mathbb{R}^{+})}ds+C\int_0^{t}\|\langle \eta\rangle^{6+6\delta}v^{T}_{\eta} G\|_{L_{x}^1(\mathbb{R}^{+})}ds,
\end{align}
which implies \eqref{xiuppo999iu89}, if $\varepsilon_1$ is sufficiently small.
\end{proof}

\section{Proof of Theorem \ref{mainthm}}\label{xhzuyo898}
Now we will prove Theorem \ref{mainthm} by some bootstrap argument. Assume that $u$ is a classical solution to the initial-boundary value problem \eqref{quasiwave}, \eqref{boundary}, \eqref{initial}. We will show that there exist positive constants $\varepsilon_0$ and $A$ such that
\begin{align}\label{rgty6788888}
\sup_{0\leq s\leq t}E_4(u(s))+\mathcal {E}_4(u(t))+\sup_{0\leq s\leq t}\mathbb{E}_3(u(s))+\mathscr{E}_3(u(t))\leq A^2\varepsilon^2
\end{align}
 under
the assumption
\begin{align}\label{dddddd888}
\sup_{0\leq s\leq t}E_4(u(s))+\mathcal {E}_4(u(t))+\sup_{0\leq s\leq t}\mathbb{E}_3(u(s))+\mathscr{E}_3(u(t))\leq 4A^2\varepsilon^2,
\end{align}
where $0<\varepsilon\leq \varepsilon_0$.

First, note that \eqref{xxxjkdddd900}, \eqref{xuyao890hj} and \eqref{quasiwave} imply
\begin{align}
\label{quasilinear_id} E_4(u(0))+\mathbb{E}_3(u(0)) \leq C_0 \varepsilon^2 \text{,}
\end{align}
for some constant $C_0 > 0$.
Also, Lemma \ref{xuyao8DD8899} and \eqref{dddddd888} imply \eqref{POINT12}, with $\nu_1 \simeq \varepsilon$, hence the conclusion of Lemma \ref{point222} holds throughout.

\subsection{Control of high order energy}
In order to estimate $\sup\limits_{0\leq s\leq t}{E}_4(u(s))$ and ${\mathcal {E}}_4(u(t))$,
in view of \eqref{GHUER444}, \eqref{GHUEFFFFR444}, Lemma \ref{porpkoo} and Lemma \ref{porpkoo1}, we have
\begin{align}\label{XUYAO1}
\sup_{0\leq s\leq t}E_4(u(s))\leq C\sup_{0\leq s\leq t}{E}_4^2(u(s))+C\sup_{0\leq s\leq t}\widetilde{E}_4(u(s))
\end{align}
and
\begin{align}\label{XUYAO2}
\mathcal {E}_4(u(t))\leq C\sup_{0\leq s\leq t}{E}_4(u(s)){\mathcal {E}}_4(u(t))+C\widetilde{\mathcal {E}}_4(u(t)),
\end{align}
thus our task is to control $\sup\limits_{0\leq s\leq t}\widetilde{E}_4(u(s))$ and $\widetilde{\mathcal {E}}_4(u(t))$.

From the system \eqref{quasiwave}, for $k\leq 3$
we have
\begin{align}\label{fgt}
\partial_t^ku_{\xi\eta}&=A_1\partial_t^ku_{\xi\eta}+A_2\partial_t^ku_{\xi\xi}+A_3\partial_t^ku_{\eta\eta}+G_k,
\end{align}
where
\begin{align}
G_k&=\sum_{l=1}^{k}\lambda_{kl}\big(\partial^{l}_tA_1\partial^{k-l}_tu_{\xi\eta}+\partial^{l}_tA_2\partial^{k-l}_tu_{\xi\xi}
+\partial^{l}_tA_3\partial^{k-l}_tu_{\eta\eta}\big)+\partial^{k}_tF,
\end{align}
$\lambda_{kl}$ are some constants.

In view of the homogeneous boundary condition \eqref{boundary}, we have $\partial_t^ku(t,0)=0$. Then
it follows from \eqref{fgt} and Lemma \ref{xuoaj99hyyty} that
\begin{align}\label{xiuppopppppiu8}
&\sup_{0\leq s\leq t}\widetilde{{E}}_4(u(s))+ {\widetilde{\mathcal{E}}_4}(u(t))\nonumber\\
&\leq C E_4(u(0))+C\sum_{k=0}^{3}\int_0^{t}\|\langle \xi\rangle^{2+2\delta}\partial_t^ku^{T}_{\xi}G_k\|_{L_{x}^1(\mathbb{R}^{+})}ds\nonumber\\
&~~~~~~~~~~~~~~~~~~+C\sum_{k=0}^{3}\int_0^{t}\|\langle \eta\rangle^{2+2\delta}\partial_t^{k}u^{T}_{\eta} G_k\|_{L_{x}^1(\mathbb{R}^{+})}ds.
\end{align}
For $k\leq 3$, Lemma \ref{point222} implies
\begin{align}\label{rfghy7888}
|G_k|\leq C\sum_{|b|\leq 2}|Z^{b}u_{\xi}| \sum_{|c|\leq 3}|Z^{c}u_{\eta}|+C\sum_{|b|\leq 3}|Z^{b}u_{\xi}| \sum_{|c|\leq 2}|Z^{c}u_{\eta}|.
\end{align}
It also holds that
\begin{align}\label{xiuppoppppttpiu8}
\int_0^{t}\|\langle \xi\rangle^{2+2\delta}\partial_t^ku^{T}_{\xi}G_k\|_{L_{x}^1(\mathbb{R}^{+})}ds
&\leq \|\langle \eta\rangle^{-\frac{1+\delta}{2}}\langle \xi\rangle^{1+\delta}\partial_t^ku_{\xi}\|_{L^2_{s,x}}\|\langle \eta\rangle^{\frac{1+\delta}{2}}\langle \xi\rangle^{1+\delta}G_k\|_{L^2_{s,x}}\nonumber\\
&\leq C\mathcal {E}_4^{1/2}(u(t))\|\langle \eta\rangle^{\frac{1+\delta}{2}}\langle \xi\rangle^{1+\delta}G_k\|_{L^2_{s,x}}.
\end{align}
Then by \eqref{rfghy7888} and Lemma \ref{xuyao8DD8899} we obtain
\begin{align}\label{xyu9dddd00}
&\|\langle \eta\rangle^{\frac{1+\delta}{2}}\langle \xi\rangle^{1+\delta}G_k\|_{L^2_{s,x}}\nonumber\\
&\leq C
\sum_{|b|\leq 2}\|\langle \eta\rangle^{-\frac{1+\delta}{2}}\langle \xi\rangle^{1+\delta}Z^{b}u_{\xi}\|_{L^2_{s}L^{\infty}_{x}}  \sum_{|c|\leq 3}\|\langle\eta\rangle^{1+\delta}Z^{c}u_{\eta}\|_{L^{\infty}_{s}L^2_{x}}\nonumber\\
&+C
\sum_{|b|\leq 3}\|\langle \eta\rangle^{-\frac{1+\delta}{2}}\langle \xi\rangle^{1+\delta}Z^{b}u_{\xi}\|_{L^2_{s,x}}\sum_{|c|\leq 2}\|\langle\eta\rangle^{1+\delta}Z^{c}u_{\eta}\|_{L^{\infty}_{s,x}}\nonumber\\
&\leq C\sup_{0\leq s\leq t}E_4^{1/2}(u(s))\mathcal {E}_4^{1/2}(u(t)).
\end{align}
It follows from \eqref{xiuppoppppttpiu8} and \eqref{xyu9dddd00} that
\begin{align}\label{xiuppoppppjjjjttpiu8}
\int_0^{t}\|\langle \xi\rangle^{2+2\delta}\partial_t^ku^{T}_{\xi}G_k\|_{L_{x}^1(\mathbb{R}^{+})}ds
\leq C\sup_{0\leq s\leq t}E_4^{1/2}(u(s))\mathcal {E}_4(u(t)).
\end{align}
Similarly, we can also show
\begin{align}\label{xiuppopppdddpjjjjttpiu8}
\int_0^{t}\|\langle \eta\rangle^{2+2\delta}\partial_t^ku^{T}_{\eta}G_k\|_{L_{x}^1(\mathbb{R}^{+})}ds
\leq C\sup_{0\leq s\leq t}E_4^{1/2}(u(s))\mathcal {E}_4(u(t)).
\end{align}
Thus the combination of \eqref{xiuppopppppiu8}, \eqref{xiuppoppppjjjjttpiu8} and \eqref{xiuppopppdddpjjjjttpiu8} gives
 \begin{align}\label{XUYAO3}
\sup_{0\leq s\leq t}\widetilde{{E}}_4(u(s))+ {\widetilde{\mathcal{E}}_4}(u(t))\leq C E_4(u(0))+C\sup_{0\leq s\leq t}E_4^{1/2}(u(s))\mathcal {E}_4(u(t)).
\end{align}

Finally,
in view of \eqref{XUYAO1}, \eqref{XUYAO2} and \eqref{XUYAO3}, we have
\begin{align}\label{XUYAO4}
&\sup_{0\leq s\leq t}{{E}}_4(u(s))+ {{\mathcal{E}}}_4(u(t))\nonumber\\
&\leq C E_4(u(0))+C\big(\sup_{0\leq s\leq t}E_4^{1/2}(u(s))+\sup_{0\leq s\leq t}E_4(u(s))\big)\mathcal {E}_4(u(t))+C\sup_{0\leq s\leq t}E_4^{2}(u(s)).
\end{align}

\subsection{Control of low order energy}
The remaining task is to estimate $\sup\limits_{0\leq s\leq t}\mathbb{E}_3(u(s))$ and $\mathscr{E}_3(u(t))$. By \eqref{GHUER444}, \eqref{GHUEFFFFR444}, Lemma \ref{porpkoo} and Lemma \ref{porpkoo1}, we see
\begin{align}\label{XUYAO19}
\sup_{0\leq s\leq t}\mathbb{E}_3(u(s))\leq C\sup_{0\leq s\leq t}{E}_3(u(s))\sup_{0\leq s\leq t}\mathbb{E}_3(u(s))+C\sup_{0\leq s\leq t}\widetilde{\mathbb{E}}_3(u(s))
\end{align}
and
\begin{align}\label{XUYAO29}
\mathscr{E}_3(u(t))\leq C\sup_{0\leq s\leq t}{E}_3(u(s)){\mathscr{E}}_3(u(t))+C\widetilde{\mathscr{E}}_3(u(t)).
\end{align}
Thus we need to control $\sup\limits_{0\leq s\leq t}\widetilde{\mathbb{E}}_3(u(s))$ and $\widetilde{\mathscr{E}}_3(u(t))$.


Note that for $k\leq 2, ~$$\partial^{k}_tu$ satisfies the system \eqref{fgt} and the homogeneous boundary condition $\partial_t^ku(t,0)=0$. By Lemma \ref{xuoaj99hyyty9} we have
\begin{align}\label{xiuppopppppiu89}
&\sup_{0\leq s\leq t}\widetilde{{\mathbb{E}}}_3(u(s))+ {\widetilde{\mathscr{E}}_3}(u(t))\nonumber\\
&\leq C \mathbb{E}_3(u(0))+C\mathscr{E}_2^{1/2}(u(t))\mathscr{E}_3^{1/2}(u(t))\sup_{0\leq s\leq t}E^{1/2}_4(u(t))\nonumber\\
&+C\sum_{k=0}^{2}\int_0^{t}\|\langle \xi\rangle^{6+26\delta}\partial_t^ku^{T}_{\xi}G_k\|_{L_{x}^1(\mathbb{R}^{+})}ds+C\sum_{k=0}^{2}\int_0^{t}\|\langle \eta\rangle^{6+6\delta}\partial_t^{k}u^{T}_{\eta} G_k\|_{L_{x}^1(\mathbb{R}^{+})}ds.
\end{align}
For $k\leq 2$, Lemma \ref{point222} implies
\begin{align}\label{rfghy78889}
|G_k|\leq C\sum_{|b|\leq 1}|Z^{b}u_{\xi}| \sum_{|c|\leq 2}|Z^{c}u_{\eta}|+C\sum_{|b|\leq 2}|Z^{b}u_{\xi}| \sum_{|c|\leq 1}|Z^{c}u_{\eta}|.
\end{align}
We also have
\begin{align}\label{xiuppoppppttpiu89}
\int_0^{t}\|\langle \eta\rangle^{6+6\delta}\partial_t^ku^{T}_{\eta}G_k\|_{L_{x}^1(\mathbb{R}^{+})}ds
&\leq \|\langle \xi\rangle^{-\frac{1+\delta}{2}}\langle \eta\rangle^{3+3\delta}\partial_t^ku_{\eta}\|_{L^2_{s,x}}\|\langle \xi\rangle^{\frac{1+\delta}{2}}\langle \eta\rangle^{3+3\delta}G_k\|_{L^2_{s,x}}\nonumber\\
&\leq C\mathscr{E}_3^{1/2}(u(t))\|\langle \xi\rangle^{\frac{1+\delta}{2}}\langle \eta\rangle^{3+3\delta}G_k\|_{L^2_{s,x}}.
\end{align}
Then by \eqref{rfghy78889} and Lemma \ref{xuyao8DD8899} we can get
\begin{align}\label{xyu9dddd009}
&\|\langle \xi\rangle^{\frac{1+\delta}{2}}\langle \eta\rangle^{3+3\delta}G_k\|_{L^2_{s,x}}\nonumber\\
&\leq C
\sum_{|b|\leq 1}\|\langle \xi\rangle^{-\frac{1+\delta}{2}}\langle \eta\rangle^{3+3\delta}Z^{b}u_{\eta}\|_{L^2_{s}L^{\infty}_{x}}  \sum_{|c|\leq 2}\|\langle\xi\rangle^{1+\delta}Z^{c}u_{\xi}\|_{L^{\infty}_{s}L^2_{x}}\nonumber\\
&+C
\sum_{|b|\leq 2}\|\langle \xi\rangle^{-\frac{1+\delta}{2}}\langle \eta\rangle^{3+3\delta}Z^{b}u_{\eta}\|_{L^2_{s,x}}\sum_{|c|\leq 1}\|\langle\xi\rangle^{1+\delta}Z^{c}u_{\xi}\|_{L^{\infty}_{s,x}}\nonumber\\
&\leq C\sup_{0\leq s\leq t}{E}_3^{1/2}(u(s))\mathscr{E}_3^{1/2}(u(t)).
\end{align}
It follows from \eqref{xiuppoppppttpiu89} and \eqref{xyu9dddd009} that
\begin{align}\label{xiuppoppppjjjjttpiu89}
\int_0^{t}\|\langle \eta\rangle^{6+6\delta}\partial_t^ku^{T}_{\eta}G_k\|_{L_{x}^1(\mathbb{R}^{+})}ds
\leq C\sup_{0\leq s\leq t}{E}_3^{1/2}(u(s))\mathscr{E}_3(u(t)).
\end{align}
Similarly, we can also show
\begin{align}\label{xiuppopppdddpjjjjttpiu89}
\int_0^{t}\|\langle \xi\rangle^{6+6\delta}\partial_t^ku^{T}_{\xi}G_k\|_{L_{x}^1(\mathbb{R}^{+})}ds
\leq C\sup_{0\leq s\leq t}{E}_3^{1/2}(u(s))\mathscr{E}_3(u(t)).
\end{align}
Thus the combination of \eqref{xiuppopppppiu89}, \eqref{xiuppoppppjjjjttpiu89} and \eqref{xiuppopppdddpjjjjttpiu89} gives
 \begin{align}\label{XUYAO39}
&\sup_{0\leq s\leq t}\widetilde{{\mathbb{E}}}_3(u(s))+ {\widetilde{\mathscr{E}}_3}(u(t))\leq C \mathbb{E}_3(u(0))+C\sup_{0\leq s\leq t}{E}_4^{1/2}(u(s))\mathscr{E}_3(u(t)).
\end{align}

Finally,
by \eqref{XUYAO19}, \eqref{XUYAO29} and \eqref{XUYAO39}, we have
\begin{align}\label{XUYAO49}
&\sup_{0\leq s\leq t}\mathbb{E}_3(u(s))+\mathscr{E}_3(u(t))\nonumber\\
&\leq C\mathbb{E}_3(u(0))+C\big(\sup_{0\leq s\leq t}E_4^{1/2}(u(s))+\sup_{0\leq s\leq t}E_3(u(s))\big)\mathscr{E}_3(u(t))\nonumber\\
&+C\sup_{0\leq s\leq t}{E}_3(u(s))\sup_{0\leq s\leq t}\mathbb{E}_3(u(s)).
\end{align}

\subsection{Global existence}
Now we have arrived at the coupled energy estimates \eqref{XUYAO4} and \eqref{XUYAO49}.
Thus under the assumption \eqref{dddddd888}, we have
\begin{align}
&\sup_{0\leq s\leq t}E_4(u(s))+\mathcal {E}_4(u(t))+\sup_{0\leq s\leq t}\mathbb{E}_3(u(s))+\mathscr{E}_3(u(t))\nonumber\\
&\leq C_1\varepsilon^2+8C_1A^3\varepsilon^3+16C_1A^4\varepsilon^4.
\end{align}
Note that \eqref{quasilinear_id} holds.
%
Taking $A^2=4\max\{C_0,{C_1}\}$ and $\varepsilon_0$ so small that
\begin{align}
32C_1A\varepsilon_0+64C_1A^2\varepsilon^2_0\leq 1,
\end{align}
for any $\varepsilon$ with $0<\varepsilon\leq \varepsilon_0$, we have
\begin{align}\label{351}
\sup_{0\leq s\leq t}E_4(u(s))+\mathcal {E}_4(u(t))+\sup_{0\leq s\leq t}\mathbb{E}_3(u(s))+\mathscr{E}_3(u(t))\leq A^2\varepsilon^2,
\end{align}
which completes the proof of Theorem \ref{mainthm}.

\section*{Acknowledgements}
The author is supported by National Natural Science Foundation of China (No.12371217) and the Fundamental Research Funds for the Central Universities (No. 2232022D-27).

\vspace{0.3cm}
\noindent{\bf {Declarations}}

\vspace{0.3cm}
\noindent
{\bf {Conflict of interest}} On behalf of all authors, the corresponding author states that there is no conflict of interest.

\vspace{0.3cm}
\noindent
{\bf {
Data Availability}} There is no data associated with this work.


\end{document}